\documentclass{amsart}
\usepackage[margin=2cm]{geometry}
\usepackage{fontenc,amsmath,mathrsfs, amsfonts, amsthm ,amssymb, MnSymbol, enumerate, multicol, hyperref, geroENG}
\usepackage[english]{babel}
\allowdisplaybreaks
\usepackage{times}
\begin{document}
\title{Metrical properties of weighted products of consecutive L\"uroth digits}
\author{Adam Brown-Sarre}
\address{Adam Brown-Sarre, Department of Mathematical and Physical Sciences,  La Trobe University, Bendigo 3552, Australia. }
\email{20356213@students.latrobe.edu.au}

\author{Gerardo Gonz\'alez Robert}
\address{Gerardo Gonz\'alez Robert, Department of Mathematical and Physical Sciences,  La Trobe University, Bendigo 3552, Australia. }
\email{G.Robert@latrobe.edu.au}

\author{Mumtaz Hussain}
\address{Mumtaz Hussain,  Department of Mathematical and Physical Sciences,  La Trobe University, Bendigo 3552, Australia. }
\email{m.hussain@latrobe.edu.au}

\date{}

\begin{abstract}
The L\"uroth expansion of a real number $x\in (0,1]$ is the series
\[
x=
\frac{1}{d_1} + \frac{1}{d_1(d_1-1)d_2} + \frac{1}{d_1(d_1-1)d_2(d_2-1)d_3} + \cdots,
\]
with $d_j\in\mathbb{N}_{\geq 2}$ for all $j\in\mathbb{N}$. Given $m\in \mathbb{N}$, $\mathbf{t}=(t_0,\ldots, t_{m-1})\in\mathbb{R}_{>0}^{m-1}$ and any function $\Psi:\mathbb{N}\to (1,\infty)$, define
\[
\clE_{\mathbf{t}}(\Psi)\colon= \left\{ x\in (0,1]: d_n^{t_0} \cdots d_{n+m}^{t_{m-1}}\geq \Psi(n) \text{ for infinitely many} \ n \in\mathbb{N} \right\}.
\]
We establish a Lebesgue measure dichotomy statement (a zero-one law) for $\clE_{\mathbf{t}}(\Psi)$ under a natural non-removable condition $\liminf_{n\to\infty} \Psi(n)>~1$.  Let $B$ be given by 
\[
\log B \colon= \liminf_{n\to\infty} \frac{\log(\Psi(n))}{n}. 
\]
For any $m\in\mathbb{N}$, we compute the Hausdorff dimension of $\clE_{\mathbf{t}}(\Psi)$ when either $B=1$ or $B=\infty$. We also compute the Hausdorff dimension of $\clE_{\mathbf{t}}(\Psi)$ when $1<B< \infty$ for $m=2$.
\end{abstract}
\maketitle
\section{Introduction}
J. L\"uroth showed in 1883 that every $x\in (0,1]$ can be uniquely expressed as a series of the form
\begin{equation*}\label{Eq:Lu-02}
x= \frac{1}{d_1(x)} + \frac{1}{d_1(x)(d_1(x)-1)d_2(x)} + \frac{1}{d_1(x)(d_1(x)-1)d_2(x)(d_2(x)-1)d_3(x)} + \cdots,
\end{equation*}
where each $d_j(x)$ is an integer at least $2$. Analogous to regular continued fractions with the Gauss map, L\"uroth series can be interpreted through dynamical means using the L\"uroth map $\scL$ (see Section \ref{SEC:LUROTH} for its definition).

The metrical theory of L\"uroth expansions has been thoroughly studied. In \cite{JagDev1969}, H. Jager and C. de Vroedt showed that the Lebesgue measure $\lambda$ on $(0,1]$ is $\scL$-ergodic. They also noted that the digits in the L\"uroth expansion are independent when regarded as random variables. Hence, the L\"uroth series analogue of the Borel-Bernstein Theorem (for regular continued fractions) \cite[Theorem 30]{Khi1997} is an immediate consequence of the classical Borel-Cantelli Lemma (see, for example, \cite[Theorem 4.18]{Kal2021}) and the following observation:
\[
\lambda\left(\left\{ x\in (0,1]: d_n(x)\geq m\right\}\right)=\frac{1}{m}
\quad
\text{ for all } n\in\mathbb{N} \text{ and all } m\in\mathbb{N}_{\geq 2}.
\]


Throughout this paper, $\Psi:\mathbb{N}\to\mathbb{R}_{>0}$ denotes a positive function. The Borel-Bernstein Theorem for Lüroth series provides the Lebesgue measure $\lambda$ of the set
\[
\clE(\Psi)
\colon=
\left\{ x\in (0,1]: d_n(x) \geq \Psi(n) \ \text{\rm for infinitely many } n\in\mathbb{N}\right\}.
\]
\begin{teo01}[{{\cite[Theorem 2.1]{JagDev1969}}}]\label{TEO:BB:L}
The Lebesgue measure of $\clE(\Psi)$ is given by
\[
\lambda\left( \clE(\Psi)\right)
= 
\begin{cases}
0, &\text{\rm if } \quad \displaystyle\sum_{n=1}^{\infty} \Psi(n)^{-1} <\infty,\\[2ex]
1, &\text{\rm  if }\quad  \displaystyle\sum_{n=1}^{\infty} \Psi(n)^{-1} = \infty.
\end{cases}
\]
\end{teo01}
Theorem \ref{TEO:BB:L} appeared first in \cite{JagDev1969}, although $\Psi$ was unnecessarily assumed to be increasing. 

It is well-known that to Hausdorff dimension is an appropriate notion to distinguish between the null sets (Lebesgue measure zero sets). L. Shen  \cite{She2017}  calculated the Hausdorff dimension of the set  $ \clE(\Psi)$ which is the L\"uroth analogue of a celebrated theorem of B. Wang and J. Wu \cite{WanWu2008}. To analyse the Hausdorff dimension of $\clE(\Psi)$ and related sets, define 
\begin{equation}\label{EQ:DefBb}
\log B
=
\liminf_{n\to\infty} \frac{\log \Psi(n)}{n}
\;\text{ and }\;
\log b 
=
\liminf_{n\to\infty} \frac{\log\log \Psi(n)}{n}.
\end{equation}
\begin{teo01}[{\cite[Theorem 4.2]{She2017}} ]\label{TEO:WAWU:L}
The Hausdorff dimension of $\clE(\Psi)$ is as follows:
\begin{enumerate}[1.]
\item If $B=1$, then $\dimh \clE(\Psi)=1$.
\item If $1<B<\infty$, then $\dimh \clE(\Psi)=s(B)$, where $s=s(B)$ is the solution of 
\[
\sum_{k=2}^{\infty} \left(\frac{1}{B k(k-1)}\right)^s=1.
\]
\item When $B=\infty$, we have three cases:
\begin{enumerate}[\rm (i)]
\item If $b=1$, then $\dimh \clE(\Psi)=\frac{1}{2}$.
\item If $1<b<\infty$, then $\dimh \clE(\Psi)=\frac{1}{1+b}$.
\item If $b=\infty$, then $\dimh \clE(\Psi)=0$.
\end{enumerate}
\end{enumerate}
\end{teo01}

In \cite{TanZho2021}, B. Tan and Q. Zhou extended Shen's result by  considering the product of two consecutive partial quotients. Define the set
\[
\clE_{(1,1)}(\Psi)
\colon
=\left\{ x\in (0,1]: d_n(x)d_{n+1}(x)\geq \Psi(n) \ \text{\rm for infinitely many } n\in\mathbb{N}\right\}.
\]
\begin{teo01}[{\cite[Lemma 3.1]{TanZho2021}}]\label{TEO:QZ:L}
Let $t=t(B)$ be the unique solution of the equation
\[
\sum_{d=2}^{\infty} \frac{1}{B^{2s}d^s(d-1)^s} = 1.
\]
If $B$ and $b$ are given by  \eqref{EQ:DefBb}, then 
\[
\dimh\clE_{(1,1)}(\phi)
=
\begin{cases}
1, &\text{\rm if }\quad B=1, \\
t(B), &\text{\rm if }\quad 1<B<\infty, \\
\frac{1}{1+b}, &\text{\rm if }\quad  B=\infty.
\end{cases}
\]
\end{teo01}
In this paper, we consider the weighted product of consecutive partial quotients and establish the Lebesgue measure and Hausdorff dimension for the corresponding limsup set. Given $m\in\mathbb{N}$ and $\mathbf{t}=(t_0,\ldots, t_{m-1})\in\mathbb{R}_{>0}^m$, define the set
\[
\clE_{\mathbf{t}}(\Psi)
\colon=
\left\{ x\in (0,1]: \prod_{i=0}^{m-1} d_{n+i}^{t_i}(x) \geq \Psi(n) \text{ for infinitely many }n\in\mathbb{N}\right\},
\]
and the numbers
\[
t \colon= \min\{t_0,t_1,\ldots, t_{m-1}\}, 
\quad
T \colon=\max\{t_0,t_1,\ldots, t_{m-1}\},  
\quad
\ell(\mathbf{t})\colon=\# \{j\in \{0,\ldots, m\}: t_j=T\}. 
\]
\begin{teo01}\label{TEO:01}
Let $m\in\mathbb{N}$ and $\mathbf{t}\in \mathbb{R}_{>0}^m$ be arbitrary. If
\begin{equation}\label{EQ:TEO:01:liminf}
\liminf_{n\to\infty} \Psi(n)>1, 
\end{equation}
then, for $B$ and $b$ given by  \eqref{EQ:DefBb}, we have
\begin{equation}\label{EQ:TEO:01:Conclusion}
\lambda\left(\clE_{\mathbf{t}}(\Psi)\right)
=
\begin{cases}
0, &\text{\rm if }\quad \displaystyle\sum_{n=1}^{\infty} \frac{\left( \log\Psi(n)\right)^{\ell(\mathbf{t}) - 1} }{\Psi(n)^{\frac{1}{T}}} < \infty, \\[2ex]
1, &\text{\rm if }\quad \displaystyle\sum_{n=1}^{\infty} \frac{\left( \log\Psi(n)\right)^{\ell(\mathbf{t}) - 1}}{\Psi(n)^{\frac{1}{T}}} = \infty.
\end{cases}
\end{equation}
\end{teo01}
Assumption \eqref{EQ:TEO:01:liminf} essentially says that there exists some $n_0\in\mathbb{R}_{>1}$ such that $\Psi(n)\geq n_0$ for all $n\in\mathbb{N}$. Theorem \ref{TEO:01} might fail without this condition (see Section \ref{SEC:PROOF:TEO:01} for the justification). 

\begin{teo01}\label{TEO:02}
Take any $m\in\mathbb{N}$ and $\mathbf{t}\in \mathbb{R}_{>0}^{m}$ and let $B$ and $b$ be given by  \eqref{EQ:DefBb}. Then
\[
\dimh\clE_{\mathbf{t}}(\Psi)
=
\begin{cases}
1, &\text{\rm if }\quad  B=1, \\
\frac{1}{b+1}, &\text{\rm if }\quad  B=\infty.
\end{cases}
\]
\end{teo01}
If $m=2$, we can compute the $\dimh \clE_{\mathbf{t}}(\Psi)$ when $1<B<\infty$.
\begin{teo01}\label{TEO:03} 
Let $B$ and $b$ be given by  \eqref{EQ:DefBb} and assume that $1<B<\infty$. For a given $\mathbf{t}=(t_0,t_1)\in\mathbb{R}_{>0}^2$, define
\[
f_{t_0,t_1}(s)\colon= \frac{s^2}{t_0t_1\max \left\{ \frac{s}{t_1} + \frac{1-s}{t_0},\frac{s}{t_0}\right\}}.
\]
The Hausdorff dimension of $\clE_{\mathbf{t}}(\Psi)$ is the unique solution of
\[
\sum_{d=2}^{\infty} \frac{1}{d^s(d-1)^s B^{f_{t_0,t_1}(s)}}=1.
\]
\end{teo01}

For completeness, it is worth mentioning that there have been abundance of work regarding the metrical properties of product of consecutive partial quotients. The results of the paper  provide full L\"uroth analogues of very recent work of Bakhtawar-Hussain-Kleinbock-Wang \cite{BakHusKleWan2022}. Note that paper \cite{BakHusKleWan2022} was a generalisation of previous works \cite{BHS2023,  HuaWuXu2020, HKWW2018, KleWad2018, WanWu2008}.

\medskip

The organization of the paper is as follows. In Section \ref{SEC:LUROTH}, we recall some basic facts on L\"uroth series. In Section \ref{SEC:PROOF:TEO:01}, we prove Theorem \ref{TEO:01}. Section \ref{SEC:PROOF:TEO:02} is dedicated to the proof of Theorem \ref{TEO:02}. In Section \ref{SEC:PROOF:TEO:03}, we prove Theorem \ref{TEO:03}. Lastly, in Section \ref{SEC:FINALREMARKS}, we give some final remarks and a conjecture.

\paragraph{Notation.} We adopt the Vinogradov symbol $\ll$ for asymptotic behavior. If $\saxu$ and $\sayu$ are two sequences of positive real numbers, we write $x_n\ll y_n$ if there exists a constant $C>0$ such that $x_n\leq C y_n$ holds for all $n\in\mathbb{N}$. When the constant $C$ depends on some parameter $m$, we write $x_n\ll_m y_n$. If we have $x_n\ll y_n$ and $y_n\ll x_n$, we write $x_n\asymp y_n$. If the implied constants depend of some parameter $m$, we write $x_n\asymp_m y_n$. We write $\scD\colon=\mathbb{N}_{\geq 2}$. If $\bfa=(a_1,\ldots, a_n)\in \scD^n$ and $\bfb=(b_1,\ldots, b_m)\in \scD^m$, then $\bfa\bfb\in\scD^{n+m}$ is $\bfa\bfb\colon=(a_1,\ldots,a_n,b_1,\ldots,b_{m})$. We denote by $\lambda$ the Lebesgue measure on $\mathbb{R}$.

\medskip

\noindent{\bf Acknowledgements} The research of Mumtaz Hussain and Gerardo González Robert is supported by the Australian Research Council Discovery Project (200100994).

\section{Elements of L\"uroth series}\label{SEC:LUROTH}

Let $d_1:(0,1]\to\scD\colon=\mathbb{N}_{\geq 2}$ be the function associating to each $x\in (0,1]$ the natural number $d_1(x)\geq 2$ determined by
\[
\frac{1}{d_1(x)}< x \leq \frac{1}{d_1(x)-1}.
\]
That is, if $\lfloor\,\cdot\,\rfloor$ represents the floor function, we define $d_1(x)\colon= \lfloor \frac{1}{x}\rfloor + 1$. The \textbf{L\"uroth map} is the function $\scL(x):[0,1]\to [0,1]$ given by
\[
\scL (x)= 
\begin{cases}
d_1(x)(d_1(x)-1)x - (d_1(x)-1), \text{ if } x\in (0,1],\\
0, \text{ if } x=0.
\end{cases}
\]
For any $x\in (0, 1]$ and $n\geq 2$, we define $d_n(x)\colon =d_1(\scL^{n-1}(x))$, the exponent denotes iteration. For any $\mathbf{d}:=(d_1,\ldots, d_n)\in\scD^{n}$, the \textbf{cylinder} of level $n$ based at $\mathbf{d}$ is the set \[
I_n(\mathbf{d})
\colon=
\left\{ x\in (0,1]: d_1(x)=d_1, \ldots, d_n(x)=d_n\right\}.
\]
The cylinders are intervals of the form\footnote{Some authors, for example  \cite{DajKra2002}, define the Lüroth map in a slightly different manner. For any $x\in (0,1)$ they choose $n\in \mathbb{N}$ such that $ n^{-1}\leq x < (n-1)^{-1}$, they map $x$ to $n(n-1)x - (n-1)$, and they leave $0$ as a fixed point. Under this definition, cylinders are intervals of the for $[\alpha,\beta)$. Our definition follows \cite{Gall976, JagDev1969} among others. The difference between these two approaches is irrelevant for our purposes, because it only affects the countable set $\QU\cap [0,1]$.} $(\alpha, \beta]$ for $0<\alpha<\beta\leq 1$. 
 
L\"uroth series induce a continuous map $\Lambda:\scD^{\mathbb{N}} \to (0,1]$ given by 
\begin{align*}
\Lambda(d_1,d_2,d_3,\ldots)
&\colon=
\langle d_1,d_2,d_3,\ldots \rangle \\
&\colon=
\frac{1}{d_1} + \frac{1}{d_1(d_1 - 1)d_2}  + \frac{1}{d_1(d_1 - 1)d_2(d_2-1)d_3} + \cdots.  
\end{align*}
Denote by $\sigma$ the left shift on $\scD^{\mathbb{N}}$. Then, the dynamical system $(\scD^{\mathbb{N}},\sigma)$ is an extension of $((0,1],\scL)$ in the sense that $\Lambda:\scD^{\mathbb{N}}\to (0,1]$ is a continuous onto map satisfying $\Lambda \circ \sigma= \scL\circ \Lambda$. 
Clearly, these systems cannot be topologically conjugated, because $\scD^{\mathbb{N}}$ is totally disconnected and $(0,1]$ is connected (see \cite[Proposition 2.1]{AlcGero2022}).

For each $(d_1,\ldots, d_n)\in\scD^n$, write
\[
\langle d_1, \ldots, d_n\rangle
\colon=
\frac{1}{d_1} + \frac{1}{d_1(d_1 - 1)d_2} + \cdots + \frac{1}{d_1(d_1 - 1)d_2(d_2-1)\cdots d_{n-1}(d_{n-1}-1)d_n}.
\]
\begin{propo01}\label{PROP:DIAMETER_INTERVALS}
For every $n\in\mathbb{N}$ and every $\bfc=(c_1,\ldots, c_n)\in \scD^n$, we have
\[
I_n(\bfc)
=
\left( \langle c_1,\ldots, c_n\rangle, \langle c_1,\ldots, c_n-1\rangle\right);
\]
and therefore,
\[
|I_n(\bfc)|
=
\prod_{j=1}^n \frac{1}{c_j(c_j-1)}.
\]
\end{propo01}
\begin{proof}
The proof follows from applying mathematical induction on $n$.
\end{proof}
\section{Proof of Theorem \ref{TEO:01} }\label{SEC:PROOF:TEO:01}
In 1967, W. Philipp published a quantitative version of the Borel-Cantelli Lemma \cite[Theorem 3]{Phi1967}. As noted by D. Kleinbock and N. Wadleigh \cite[Remark 3.2]{KleWad2018}, Philipp's theorem can be strengthened to obtain Lemma \ref{TEO:PhilipBC}. In this lemma and afterwards, we denote the number of elements in a given set $Y$ by $\#Y$.

\begin{lem01}\label{TEO:PhilipBC}
Let $(X,\scB,\mu)$ be a probability space and let $(E_n)_{n\geq 1}$ be a sequence of measurable sets. For each $N\in\mathbb{N}$ and each $t\in X$, define
\[
A(N,t)\colon =\#\left\{ n\in\{1,\ldots,N\} \colon t\in E_n\right\}
\]
and 
\[
\vphi(N)\colon = \sum_{n=1}^N \mu(E_n).
\]
Suppose that there is a summable sequence of non-negative real numbers $(C_j)_{j\geq 1}$ such that for any $k, m,n\in\mathbb{N}$ satisfying $n+k<m$ we have
\[
\mu(E_n\cap E_m)\leq \mu(E_n)\mu(E_m)+\mu(E_m)C_{m-n}.
\]
Then, $\mu$-almost every $t\in X$ satisfies
\[
A(N,t)
= 
\vphi(N) + \mathscr{O}\left(\sqrt{\vphi(N)} \log^{\frac{3}{2}  + \veps}\vphi(N) \right)
\quad\text{ for all } \veps>0.
\]
\end{lem01}

\begin{lem01}\label{LEMMA:BHKW:3.1}
Fix $k\in\mathbb{N}$. Let $(A_n)_{n\geq 1}$ be a sequence of at most countable unions of cylinders of level $k$. We have
\[
\lambda\left( \limsup_{n\to\infty} \scL^{-n}[A_n]\right)
=
\begin{cases}
0, &\text{\rm if }\quad \displaystyle\sum_{n=1}^{\infty} \lambda(A_n)< \infty, \\
1, &\text{\rm if }\quad\displaystyle\sum_{n=1}^{\infty} \lambda(A_n)=  \infty.
\end{cases}
\]
\end{lem01}
\begin{proof}
The convergence case follows from the convergence part of the Borel-Cantelli Lemma. Assume that 
\begin{equation}\label{EQ:LEMMA:BHKW:3.1:0.5}
\sum_{n=1}^{\infty} \lambda\left( A_n\right)= \infty.
\end{equation}
The $\scL$-invariance of $\lambda$ and \eqref{EQ:LEMMA:BHKW:3.1:0.5} imply 
\[
\sum_{n=1}^{\infty} \lambda\left( \scL^{-n}[A] \right)= \infty.
\]
Let $\scB$ be the Borel $\sigma$-algebra of $(0,1]$. It is shown in \cite[Equation (3.4)]{JagDev1969} that any cylinder $I$ of level $k$ satisfies
\begin{equation}\label{EQ:LEMMA:BHKW:3.1:01}
\lambda\left( \scL^{-u}[B]\cap I\right)
=
\lambda(B)\lambda(I)
\;\text{ for all } B\in\scB \text{ and } u\in\mathbb{N}_{\geq k}.
\end{equation}
Pick an arbitrary $n\in\mathbb{N}$ and write $A_n$ as an at most countable and disjoint union of cylinders of level $k$, say
\[
A_n
=
\bigcup_{j }  I_{j}^{n}.
\]
By \eqref{EQ:LEMMA:BHKW:3.1:01} and the $\scL$-invariance of $\lambda$, for any $m\in\mathbb{N}_{> n+k}$ we have 
\begin{align*}
\lambda\left( \scL^{-n}[A_n]\cap \scL^{-m}[A_m]\right)
&= 
\lambda\left( A_n \cap \scL^{-(m-n)}[A_m]\right) \\
&=
\lambda\left(  \bigcup_{j } I_j^n \cap \scL^{-(m-n)}[A_m]\right) \\
&=
\sum_{j} \lambda\left(  I_j^n \cap \scL^{-(m-n)}[A_m]\right) \\
&=
\sum_{j=1}^{\infty} \lambda\left( I_j^n\right) \lambda(A_m) =\lambda(A_n)\lambda(A_m).
\end{align*}
The result now follows from Lemma \ref{TEO:PhilipBC}.
\end{proof}
We start with an estimate in the spirit of A. Khinchin's proof of the existence of the Lévy-Khinchin constant \cite[Theorem 31]{Khi1997}.
\begin{lem01}\label{LE:EstimLmm}
If $m\in\mathbb{N}$ and $g\geq 2^m$, then
\[
\idotsint\limits_{\substack{x_1\cdots x_m> g \\ x_1,\ldots, x_m\geq 2}} \frac{\ud x_1\cdots \ud x_m}{x_1^2\cdots x_m^2}
\asymp_{m}
\frac{\log^{m-1} g}{g}.
\]
\end{lem01}
\begin{proof}
The proof is by induction on $m$. The base $m=1$ is clear: if $g>2$, then
\[
\int\limits_{x_1>g}\frac{\ud x_1}{x_1^2}
=
\frac{1}{g}.
\]
Assume that for $m=M\in \mathbb{N}$, we have
\[
\idotsint\limits_{\substack{x_1\cdots x_M> g \\ x_1,\ldots, x_M\geq 2}} \frac{\ud x_1\cdots \ud x_M}{x_1^2\cdots x_M^2}
\asymp_M
\frac{\log^{M-1}g}{g}
\;\text{ for all } g>2^{M}.
\]
Let $g> 2^{M+1}$ be a real number. Then,
\begin{align}
&\idotsint\limits_{\substack{x_1\cdots x_Mx_{M+1}> g \\ x_1,\ldots, x_M,  x_{M+1} \geq 2}}  \frac{\ud x_1\cdots  \ud x_M \ud x_{M+1}}{x_1^2\cdots x_M^2x_{M+1}^2} = \nonumber \\
&=
\int_{2}^{\frac{g}{2^M}} \frac{1}{x_{M+1}^2}\idotsint\limits_{\substack{x_1\cdots x_M> \frac{g}{x_{M+1}} \\ x_1,\ldots, x_M\geq 2}} \frac{\ud x_1\cdots \ud x_M}{x_1^2\cdots x_M^2} \ud x_{M+1}
+
\int_{\frac{g}{2^M}}^{\infty} \frac{1}{x_{M+1}^2}\idotsint\limits_{\substack{x_1\cdots x_M> \frac{g}{x_{M+1}} \\ x_1,\ldots, x_M\geq 2}} \frac{\ud x_1\cdots \ud x_M}{x_1^2\cdots x_M^2} \ud x_{M+1}. \label{EQ:LEM:KHIN:01}
\end{align}
When $\frac{g}{2^M}<x_{M+1}$ and $\min\{x_1,\ldots, x_M\}\geq 2$, we have $x_1\cdots x_M x_{M+1} \geq 2^M > \frac{g}{x_{M+1}}$, so the second term in \eqref{EQ:LEM:KHIN:01} is 
\begin{equation}\label{EQ:PF:TEO01:01}
\int_{\frac{g}{2^M}}^{\infty} \frac{1}{x_{M+1}^2}\idotsint\limits_{\substack{x_1\cdots x_M> \frac{g}{x_{M+1}} \\ x_1,\ldots, x_M\geq 2}} \frac{\ud x_1\cdots \ud x_M}{x_1^2\cdots x_M^2} \ud x_{M+1}
=
\frac{1}{2^M}
\int_{\frac{g}{2^M}}^{\infty} \frac{1}{x_{M+1}^2} \ud x_{M+1}
=
\frac{1}{g}.
\end{equation}
If $2\leq x_{M+1} \leq \frac{g}{2^M}$, then $2^M< \frac{g}{x_{M+1}}$ and the induction hypothesis leads us to
\begin{align*}
\int_{2}^{\frac{g}{2^M}} \frac{1}{x_{M+1}^2}\idotsint\limits_{\substack{x_1\cdots x_M> \frac{g}{x_{M+1}} \\ x_1,\ldots, x_M\geq 2}} \frac{\ud x_1\cdots \ud x_M}{x_1^2\cdots x_M^2} \ud x_{M+1}
&\asymp_M 
\int_{2}^{\frac{g}{2^{M}}} \frac{1}{x_{M+1}^2} \frac{\log^{M-1}\left(\frac{g}{x_{M+1}} \right)}{\left(\frac{g}{x_{M+1}} \right)} \ud x_{M+1} \\
&=
\frac{1}{g}
\int_{2}^{\frac{g}{2^{M}}}
\frac{(\log g - \log x_{M+1})^{M-1}}{x_{M+1}} \ud x_{M+1} \\
&=
\frac{\log^{M-1} g}{g}
\int_{2}^{\frac{g}{2^{M}}} \frac{1}{x_{M+1}}
\left( 1 - \frac{\log x_{M+1}}{\log g}\right)^{M-1} \ud x_{M+1}.
\end{align*}
The inequality $0< 1 - \frac{\log x_{M+1}}{\log g}<1$ holds when $2\leq x_{M+1} \leq \frac{g}{2^M}$, then
\[
\frac{\log^{M-1} g}{g}
\int_{2}^{\frac{g}{2^M}} \frac{1}{x_{M+1}^2}\idotsint\limits_{\substack{x_1\cdots x_M> \frac{g}{x_{M+1}} \\ x_1,\ldots, x_M\geq 2}} \frac{\ud x_1\cdots \ud x_M}{x_1^2\cdots x_M^2} \ud x_{M+1}
\ll _{M+1}
\frac{\log^{M} g}{\log g}.
\]
Observe that
\begin{equation}\label{EQ:LEM:KHIN:02}
2\leq x_{M+1}\leq g^{\frac{1}{2}}
\quad\text{ implies }\quad
\frac{1}{2}
\leq
1 - \frac{\log x_{M+1}}{\log g}.
\end{equation}
We consider two cases for the lower bound: $2^{2M}< g$ and $2^{2M}\geq g$. If $2^{2M}\leq g$, then $g^{\frac{1}{2}} \leq \frac{g}{2^M}$, so
\begin{align*}
\frac{1}{g}
\int_{2}^{\frac{g}{2^{M} }}
\frac{(\log g - \log x_{M+1})^{M-1}}{x_{M+1}} \ud x_{M+1} 
&\geq 
\frac{\log^{M-1} g}{g}
\int_{2}^{ \sqrt{g} }
\frac{1}{x_{M+1}} \left( 1 - \frac{\log x_{M+1}}{\log g}\right)^{M-1} \ud x_{M+1} \\
&\geq
\frac{1}{2^{M-1}}\left(\frac{1}{2}\log g - \log 2\right) 
\frac{\log^{M-1} g}{g} \\
&=
\frac{1}{2^{M-1}} \left( \frac{1}{2} - \frac{\log 2}{\log g}\right)
\frac{\log^{M} g}{g} \\
&\geq 
\frac{1}{2^{M}}\left(1 - \frac{1}{M}\right) \frac{\log^Mg}{g}.
\end{align*}
We have used $2^{2M}\leq g$ in the last inequality.

For the second case, observe that $2^{M+1}<g \leq 2^{2M}$ yields $\frac{g}{2^M}\leq g^{\frac{1}{2}}$. This means that $2\leq x_{M+1}\leq \frac{g}{2^M}$ implies $2\leq x_{M+1}\leq g^{\frac{1}{2}}$ and, by \eqref{EQ:LEM:KHIN:02},
\begin{align*}
\frac{1}{g}
\int_{2}^{\frac{g}{2^{ M }}}
\frac{(\log g - \log x_{M+1})^M}{x_{M+1}} \ud x_{M+1} 
&\geq
\frac{1}{2^M}
\frac{\log^{M-1}g}{g} 
\int_{2}^{\frac{g}{2^{ M }}} \,
\frac{1}{x_{M+1}} \ud x_{M+1} \\
&\geq 
\frac{1}{2^{M}}
\frac{\log^{M}g}{g} 
\left( 1 - \frac{ (M+1)\log 2}{\log g} \right) >0.
\end{align*} 
These estimates along with \eqref{EQ:PF:TEO01:01} yield
\begin{align*}
\idotsint\limits_{\substack{x_1\cdots x_{M+1}> g \\ x_1,\ldots,  x_{M+1} \geq 2}} \frac{\ud x_1\cdots \ud x_{M+1}}{x_1^2\cdots x_{M+1}^2}
&\asymp_{M+1}
\frac{\log^{M} g}{g}
\left( \min\left\{ 1 - \frac{(M+1)\log 2}{\log g}, 1 - \frac{1}{M}\right\} + \frac{1}{g\log^{M} g} \right) \\
&\asymp_{M+1}
\frac{\log^{M} g}{g}.
\end{align*}
\end{proof}

\begin{lem01}\label{PROPO:LEMMA3.2:BHKW}
If $m\in\mathbb{N}$, $\mathbf{t}\in \mathbb{R}_{>0}^m$, and $g\geq 2^{m t}$, then 
\[
\sum_{\substack{d_1^{t_0}\cdots d_{m}^{t_{m-1}} \geq g \\ d_1,\ldots, d_m\geq 2}} \; \prod_{j=1}^m \frac{1}{d_j(d_j-1)} 
\asymp_{m,\mathbf{t}}
\frac{\log^{\ell-1} g }{g^{ 1/T }}.
\]
\end{lem01}
\begin{proof}
Without loss of generality, we may assume that $t=t_{m-1}\leq\ldots \leq t_1 \leq t_0 = T$, then 
\begin{equation*}
g^{ 1/t_{0} }
\leq 
\ldots
\leq
g^{ 1/t_{m-2} }
\leq 
g^{1/ t_{m-1}}.
\end{equation*}
Under this assumption on $\mathbf{t}$, the function $\ell$ has the following properties:
\begin{enumerate}[i.]
\item If $\ell(t_0,\ldots, t_{m-2}, t_{m-1})=m$, then $\ell(t_0,\ldots, t_{m-2})=m-1$.
\item If $1\leq \ell(t_0,\ldots, t_{m-2}, t_{m-1})\leq m-1$, then $\ell(t_0,\ldots, t_{m-2}, t_{m-1}) = \ell(t_0,\ldots, t_{m-2})$.
\end{enumerate}

For every $m$, $\mathbf{t}$, and $g$ as in the statement, define
\[
S_m(\mathbf{t};g)
\colon=
\sum_{\substack{d_1^{t_0}\cdots d_{m}^{t_{m-1}} \geq g \\ d_1,\ldots, d_m\geq 2}} \; \prod_{j=1}^m \frac{1}{d_j(d_j-1)}.
\]
As such, we aim to show that $S_m(\mathbf{t};g) \asymp_{m,\mathbf{t}} (\log^{\ell(\mathbf{t})-1} g)/g^{ 1/T }$.
The proof is by induction on $m$. For $m=1$ and $g>2^{t_0}$, we have 
\[
\frac{g^{1/t_0}}{2}
\leq 
\lceil g^{1/t_0}\rceil -1 
< g^{1/t_0}
\]
and so
\[
S_1(t_0;g)=
\sum_{d_1^{t_0}>g} \frac{1}{d_1(d_1-1)}
=
\sum_{d_1>\lceil g^{1/t_0}\rceil} \frac{1}{d_1(d_1-1)}
=
\frac{1}{\lceil g^{1/t_0}\rceil -1}
\asymp 
\frac{1}{g^{1/t_0}}
\]
Assume that for some $m=M\in \mathbb{N}$, every $\mathbf{t}\in\mathbb{R}_{>0}^m$ satisfies
\begin{equation}\label{Le34:IH}
S_M(\mathbf{t};g)
\asymp_{\mathbf{t},M}
\frac{\log^{\ell(\mathbf{t})-1}g}{g^{1/T}}
\;\text{ for all }g\geq 2^{mt_{M-1}}.
\end{equation}
If $\ell(\mathbf{t})=M+1$, then $S_{M+1}(\mathbf{t};g)=S_m((1,\ldots,1); g^{1/t})$ and Lemma \ref{LE:EstimLmm} gives the result. Suppose that $1\leq \ell(\mathbf{t})\leq M$. Write 
\[
u(M+1,g)\colon=
\begin{cases}
\left\lfloor \frac{g^{1/t_M}}{2^{Mt_{M-1}/t_{M}}}\right\rfloor, \text{ if } \frac{g^{1/t_M}}{2^{Mt_{M-1}/t_{M}}} \not\in\mathbb{N},\\[2ex]
\frac{g^{1/t_M}}{2^{Mt_{M-1}/t_{M}}}-1, \text{ otherwise.}
\end{cases}
\]
On the one hand, when $d_{M+1}\geq u(M+1,g)+1$, every $(d_1,\ldots, d_M)\in\scD^M$ satisfies
\[
d_1^{t_0} \cdots d_M^{t_{M-1}} 
\geq  
2^{Mt_{M-1}}
>
\frac{g}{(u(M+1,g)+1)^{t_M}}
\geq  
\frac{g}{d_{M+1}^{t_M}}.
\]
As a consequence, we may express $S_{M+1}(\mathbf{t};g)$ as follows:
\begin{align}
S_{M+1}(\mathbf{t};g)
&=
\sum_{d_{M+1}=2}^{u(M+1,g)} \frac{1}{d_{M+1}(d_{M+1}-1)} S_M\left( (t_0,\ldots,t_{M-1}), \frac{g}{d_M^{t_{M+1}}} \right) + \nonumber \\ 
&\quad+
\sum_{d_{M+1}=u(M+1,g)+1}^{\infty} \frac{1}{d_{M+1}(d_{M+1}-1)} \sum_{d_1,\ldots, d_M\geq 2} \prod_{j=1}^M \frac{1}{d_j(d_j-1)}. \label{EQ:LEMMA3.2:03}
\end{align}
We apply the induction hypothesis \eqref{Le34:IH} on the first term in \eqref{EQ:LEMMA3.2:03} and use $\ell(\mathbf{t})=\ell(t_0,\ldots, t_{M-1})$ to get
\begin{align*}
\sum_{d_{M+1}=2}^{u(M+1,g)} & \frac{1}{d_{M+1}(d_{M+1}-1)} S_M\left( (t_0,\ldots,t_{M-1}), \frac{g}{d_{M+1}^{t_{M}}} \right)\asymp_{\mathbf{t},M}\\
&\asymp_{\mathbf{t},M}
\sum_{d_{M+1}=2}^{u(M+1,g)} 
\frac{\log^{\ell(\mathbf{t})-1}(g)\left( 1 - \frac{t_M\log d_{M+1}}{\log g}\right)^{\ell(\mathbf{t})-1} }{d_{M+1}(d_{M+1}-1)\left(\frac{g}{d_{M+1}^{t_M}} \right)^{1/T}} \\
&\asymp_{\mathbf{t},M}
\frac{\log^{\ell(\mathbf{t}) - 1}(g)}{g^{1/T}}
\sum_{d_{M+1}=2}^{u(M+1,g)} \frac{1}{d_{M+1}^{2 - t_{M}/T} } 
\left( 1 - \frac{t_M\log d_{M+1}}{\log g}\right)^{\ell(\mathbf{t})-1}.
\end{align*}
Since $2- t_M/ T>1$, the last expression satisfies
\begin{align*}
\frac{\log^{\ell(\mathbf{t}) - 1}(g)}{g^{1/T}}
\sum_{d_{M+1}=2}^{u(M+1,g)} \frac{1}{d_{M+1}^{2 - t_{M}/T} } 
\left( 1 - \frac{t_M\log d_{M+1}}{\log g}\right)^{\ell(\mathbf{t}) - 1} 
&\leq 
\frac{\log^{\ell(\mathbf{t}) - 1}(g)}{g^{1/T}} 
\sum_{d_{M+1}=2}^{u(M+1,g)} \frac{1}{d_{M+1}^{2 - t_{M}/T} } \\
&\asymp_{\mathbf{t},M}
\frac{\log^{\ell(\mathbf{t}) - 1}(g)}{g^{1/T}}.
\end{align*}
We now obtain the lower estimate. If $g\leq 2^{2Mt_{M-1} + 2t_M}$, then
\[
\log g \leq (2Mt_{M-1} + 2t_M) \log 2.
\]
Therefore, since $d_{M+1} \leq u(M+1;g)$, we have
\[
\frac{t_M\log d_{M+1}}{\log g}
< 
1 - \frac{M t_{M-1}\log 2}{\log g}
\leq
1 - \frac{M t_{M-1}\log 2}{(2Mt_{M-1} + 2t_M) \log 2},
\]
which means
\[
\frac{M t_{M-1}}{2Mt_{M-1} + 2t_M}
\leq
1-\frac{t_M\log d_{M+1}}{\log g}
\]
and we conclude
\[
\frac{\log^{\ell(\mathbf{t}) - 1}(g)}{g^{1/T}}
\sum_{d_{M+1}=2}^{u(M+1,g)} \frac{1}{d_{M+1}^{2 - t_{M}/T} } 
\left( 1 - \frac{t_M\log d_{M+1}}{\log g}\right)^{\ell(\mathbf{t}) - 1} 
\gg_{\mathbf{t},M+1}
\frac{1}{4}
\frac{\log^{\ell(\mathbf{t}) - 1}(g)}{g^{1/T}}.
\]
If $g\geq 2^{2Mt_{M-1} + 2t_M}$, then $g^{1/(2t_M)}\geq 2^{1+Mt_{M-1}/t_M}$ and $2<u(M+1,g^{1/2})<u(M+1,g)$, so
\begin{align*}
\frac{\log^{\ell(\mathbf{t})-1}(g)}{g^{1/T}}
\sum_{d_{M+1}=2}^{u(M+1,g)} \frac{1}{d_{M+1}^{2 - t_{M}/T} } 
\left( 1 - \frac{t_M\log d_{M+1}}{\log g}\right)^{\ell(\mathbf{t})-1} 
&\geq
\frac{1}{2^{\ell(\mathbf{t})-1}}
\frac{\log^{\ell(\mathbf{t})-1}(g)}{g^{1/T}} 
\sum_{d_{M+1}=2}^{u(M+1,g^{1/2})} \frac{1}{d_{M+1}^{2 - t_{M}/T} }  \\
&\geq 
\frac{1}{2^{\ell(\mathbf{t})+1}}
\frac{\log^{\ell(\mathbf{t})-1}(g)}{g^{1/T}}\quad (cfr. \eqref{EQ:LEM:KHIN:02}).
\end{align*}
\end{proof}

\begin{proof}[Proof of Theorem \ref{TEO:01}]
If there are infinitely many $n\in\mathbb{N}$  such that $1< \Psi(n)\leq 2^{mt}$,  we can pick a real number $x$ and a strictly increasing sequence of natural numbers $(n_j)_{j\geq 1}$ such that 
\[
1
<
x
<\Psi(n_j)
\leq 2^{mt}
\;\text{ for all } j\in\mathbb{N}.
\]
Then, 
\[
\log^{\ell-1} x < \log^{\ell(\mathbf{t}) - 1} \Psi(n_j) 
\quad\text{ and }\quad
2^{-mt/T}< \Psi(n_j)^{-1/T}
\quad\text{ for all } j\in\mathbb{N}
\]
and, therefore,
\[
\sum_{n=1}^{\infty} \frac{\log^{\ell-1} \Psi(n) }{\Psi(n)^{1/T}}
\geq 
\sum_{j=1}^{\infty} \frac{\log^{\ell-1} \Psi(n_j) }{\Psi(n_j)^{1/T}}
=
\infty
\quad\text{ and }\quad
\clE_{\mathbf{t}}(\Psi)=(0,1].
\]
Assume that $\Psi(n)>2^{mt}$ holds for all $n\in\mathbb{N}$. For each $n\in\mathbb{N}$, define
\[
G_n^{\mathbf{t}}(\Psi)
\colon=
\left\{ \mathbf{d}\in\scD^m: d_1^{t_0}\cdots d_{m-1}^{t_m} \geq \Psi(n) \right\}
\quad\text{ and }\quad
A_n
\colon=
\bigcup_{\mathbf{d}\in G_n(\Psi)} I_m(\mathbf{d}).
\]
In view of Proposition \ref{PROP:DIAMETER_INTERVALS}, we have
\[
\lambda(A_n)
=
\sum_{\mathbf{d}\in G_n(\Psi)} \lambda(I_m(\mathbf{d}))
=
\sum_{\mathbf{d}\in G_n(\Psi)} \prod_{j=1}^m \frac{1}{d_j(d_j - 1)}
\]
and, by Lemma \ref{PROPO:LEMMA3.2:BHKW},
\begin{equation}\label{EQ:LEMMA:BHKW:3.1:02}
\lambda(A_n)
\asymp_m
\frac{ \log\Psi^{\ell(\mathbf{t})-1}(n) }{\Psi(n)^{\frac{1}{T}}}.
\end{equation}
The definitions of $\clE_{\mathbf{t}}(\Psi)$ and $(A_n)_{n=1}^{\infty}$ entail
\begin{align*}
\clE_{\mathbf{t}}(\Psi)
&=
\left\{ x\in (0,1]: \scL^{n-1}(x) \in A_n \text{ for infinitely many } n\in\mathbb{N}\right\} \\
&=
\limsup_{n\to\infty} \scL^{-n}(A_n).
\end{align*}
We deduce \eqref{EQ:TEO:01:Conclusion} from the $\scL$-invariance of $\lambda$, Lemma \ref{LEMMA:BHKW:3.1}, and \eqref{EQ:LEMMA:BHKW:3.1:02}.

Finally, we show that Theorem \ref{TEO:01} may fail if \eqref{EQ:TEO:01:liminf} is dropped. Let $t>0$ and $m\in\mathbb{N}$ be arbitrary and put $\mathbf{t}=(t,\ldots,t)\in\mathbb{R}_{>0}^m$. Choose a real number $r$ such that
\[
0 < r < \min\left\{1, \frac{1}{t}\log 2\right\}.
\]
Define $\Psi_t:\mathbb{N}\to\mathbb{R}_{\geq 1}$ by
\[
\quad
\Psi_t(n)\colon= \exp\left( r^n\right)
\quad
\text{ for all } n\in\mathbb{N}.
\]
Then, $\Psi_t(n)> 1$ for all $n\in\mathbb{N}$, $\Psi_t(n)\to 1$ as $n\to\infty$, and 
\[
\sum_{n=0}^{\infty} \frac{\log^{m-1}  \Psi_t(n) }{ \Psi_t(n)^{\frac{1}{t}} } 
=
\sum_{n=0}^{\infty} \frac{r^{n(m-1)}  }{\Psi_t(n)^{\frac{1}{t}}} 
<\infty.
\]
However, the definition of $\Psi_t$ implies $\Psi_t(n) = \exp(r^n) < 2^t$ for all $n\in\mathbb{N}$, so every $\mathbf{d}=(d_n)_{n\geq 1}\in \scD^{\mathbb{N}}$ satisfies
\[
d_n^td_{n+1}^t\cdots d_{n+m-1}^t \geq 2^{mt} > \Psi_t(n)
\quad\text{ for all } n\in\mathbb{N}.
\]
In other words, $\clE_{t}(\Psi_t)=(0,1]$ and $\lambda(\clE_{\mathbf{t}}(\Psi_t))=1$. 
\end{proof}
\section{Proof of Theorem \ref{TEO:02} }\label{SEC:PROOF:TEO:02}

We recall two results on L\"uroth series. The first one is an analogue of T. {\L}uczak's Theorem  on the Hausdorff dimension of sets of continued fractions with rapidly growing partial quotients \cite{Luc1997}. The second result is the L\"uroth analogue of a theorem by B. Wang and J. Wu \cite[Theorem 3.1]{WanWu2008}.

For every pair of real numbers $a,b$ strictly larger than $1$, define the sets
\begin{align*}
E(a,b)
&\colon=
\left\{ x \in (0,1]: d_n(x)\geq a^{b^n} \text{ for all } n\in\mathbb{N}\right\}, \\
\widetilde{E}(a,b),
&\colon=
\left\{ x \in (0,1]: d_n(x)\geq a^{b^n} \text{ for infinitely many } n\in\mathbb{N}\right\}.
\end{align*}
\begin{lem01}[{\cite[Theorem 3.1]{SheFan2011}}]\label{LEM:LUC-LUROTH}
For any $a>1$ and $b>1$, we have
\[
\dimh E(a,b) 
=
\dimh \widetilde{E}(a,b)
=
\frac{1}{b+1}.
\]
\end{lem01}
For every $B>1$, define the number
\[
s(B)
\colon=
\dimh\left\{ x\in (0,1]: d_n(x)>   B^n \text{ for infinitely many } n\in\mathbb{N}\right\}.
\]
\begin{lem01}[{\cite[Lemma 2.3]{She2017}}]\label{LEM:WWB-LUROTH}
The function $s:\mathbb{R}_{>1}\to \mathbb{R}$ is continuous, $\displaystyle\lim_{B\to\infty} s(B)=\frac{1}{2}$, and $s=s(B)$ is the only solution of
\[
\sum_{k=2}^{\infty} \left(\frac{1}{B k(k-1)}\right)^s=1.
\]
\end{lem01}
\begin{proof}[Proof of Theorem \ref{TEO:02}]
First, assume that $B=1$. From
\[
\left\{ x\in (0,1]: d_n(x)\geq \Psi(n)^{1/t_0} \text{ for infinitely many } n\in\mathbb{N}\right\}
\subseteq
\clE_{\mathbf{t}}(\Psi),
\]
$B=1$, and Theorem \ref{TEO:WAWU:L}, we conclude $\dimh\clE_{\mathbf{t}}(\Psi)=1$.

Suppose that $B=\infty$. Assume that $b>1$. Every number $c\in (1,b)$ satisfies $\frac{\log\log\Psi(n)}{n}\geq \log c$ or, equivalently, $\Psi(n) \geq e^{c^n}$ for every large $n$. Then,
\begin{align*}
\clE_{\mathbf{t}}(\Psi) 
&\subseteq
\left\{ x\in (0,1]: \prod_{i=1}^m d_{n+i}^{t_i}(x) \geq e^{c^n} \text{ for infinitely many } n\right\}\\
&\subseteq
\left\{ x\in (0,1]: d_{n+i}^{\,t_i}(x) \geq e^{\frac{c^n}{m}} \text{ for some } i\in\{1,\ldots, m\} \text{ and infinitely many } n \right\}\\
&\subseteq
\left\{ x\in (0,1]: d_{n+i}(x) \geq (e^{\frac{1}{mT}})^{c^n} \text{ for some } i\in\{1,\ldots, m\} \text{ and infinitely many } n\right\}\\
\end{align*}
and, by Lemma \ref{LEM:LUC-LUROTH},
\[
\dimh\clE_{\mathbf{t}}(\Psi)
\leq 
\frac{1}{1+c}.
\]
The previous inequality give us two implications:
\[
\text{ if } b<\infty,  \text{ then }  \dimh\clE_{\mathbf{t}}(\Psi)\leq \frac{1}{b+1},
\]
\[
\text{ if } b=\infty,  \text{ then }  \dimh\clE_{\mathbf{t}}(\Psi)=0.
\]
Now we obtain the lower bound for $\dimh \clE_{\mathbf{t}}(\Psi)$ when $1<b<\infty$. For all $c>b$, we have
\[
\left\{ x\in (0,1]: d_n(x)^{t_0}\geq e^{c^n} \text{ for all } n\in\mathbb{N}\right\}
=
\left\{ x\in (0,1]: d_n(x)\geq (e^{\frac{1}{t_0}})^{c^n} \text{ for all } n\in\mathbb{N}\right\} 
\subseteq 
\clE_{\mathbf{t}}(\Psi).
\]
Thus, applying Lemma \ref{LEM:LUC-LUROTH},
\[
\dimh \clE_{\mathbf{t}}(\Psi)
\geq 
\frac{1}{1+c} \to \frac{1}{1+b} 
\quad\text{ when }\quad c\to b.
\]
Lastly, assume that $b=1$. Then, for any $\veps>0$, we have $\Psi(n)\leq e^{(1+\veps)^n}$ infinitely often, which gives
\[
\clE_{\mathbf{t}}(\Psi)
\supseteq
\left\{ x\in (0,1]: d_{n}^{t_0}(x)\geq e^{(1+\veps)^n} \text{ for infinitely many } n\in\mathbb{N}\right\}
\]
and, by Lemma \ref{LEM:LUC-LUROTH}, 
\[
\dimh \clE_{\mathbf{t}}(\Psi)
\geq
\frac{1}{2+\veps} \to \frac{1}{2}
\quad\text{ when }\quad \veps\to 0.
\]
For the upper bound, note that $B=\infty$ implies that for each $A>0$ there is some $N(A)\in\mathbb{N}$ such that $A^n< \Psi(n)$ whenever $n\geq N(A)$. Hence, for all $x\in \clE_{\mathbf{t}}(\Psi)$ we have
\[
\prod_{i=0}^{m-1} d_{n+i}^{t_i}(x)
\geq A^n
\;\text{ for infinitely many } n\in\mathbb{N}_{\geq N(A)},
\]
and, thus, 
\begin{equation*}
d_{n+1}^{t_i} \geq A^{\frac{1}{m}}
\;\text{ for some }0\leq i\leq m-1
\;\text{ and infinitely many } n\in\mathbb{N}_{\geq N(A)}.
\end{equation*}
As a consequence,
\[
\clE_{\mathbf{t}}(\Psi)
\subseteq 
\left\{ x\in (0,1]: d_n(x) \geq (A^{\frac{1}{Tm}})^{n} \text{ for infinitely many } n\in\mathbb{N}\right\},
\]
and Lemma \ref{LEM:WWB-LUROTH} guarantees
\[
\dimh \clE_{\mathbf{t}}(\Psi)
\leq
s\left( A^{\frac{1}{Tm}} \right) \to \frac{1}{2} 
\;\text{ as }\;
A\to  \infty.
\]
\end{proof}
\section{Proof of Theorem \ref{TEO:03}}\label{SEC:PROOF:TEO:03}

We obtain Theorem \ref{TEO:03} from Theorem \ref{LEM:TEO:Aux} below and its proof. In Theorem \ref{LEM:TEO:Aux}, we compute the Hausdorff dimension of $\clE_{\mathbf{t}}(\Psi)$ for any $\mathbf{t}\in \mathbb{R}_{>0}^2$ and a particular choice of $\Psi$. For each $B\in (1,\infty)$, define the set 
\[
\clE_{\mathbf{t}}(B)
\colon=
\left\{ x\in (0,1]: d_n^{t_0}(x)d_{n+1}^{t_1}(x)\geq B^n \text{ for infinitely many } n\in\mathbb{N}\right\}.
\]
\begin{teo01}\label{LEM:TEO:Aux}
Take any $B>1$, $\mathbf{t}=(t_0,t_1)\in\mathbb{R}_{>0}^2$, and let $f_{t_0,t_1}$ be as in Theorem \ref{TEO:03}. The Hausdorff dimension of $\clE_{\mathbf{t}}(B)$ is the unique solution $s=s_0(B)$ of the equation
\[ 
\sum_{d=2}^{\infty} \frac{1}{d^s(d-1)^sB^{f_{t_{0},t_1}(s) }}=1.
\]
Moreover, the map $B\mapsto s_0(B)$ is continuous.
\end{teo01}
\subsection{Continuity of $B\mapsto s_0(B)$}
\begin{lem01}\label{Lemft0t1strictincre}
For any $(t_0,t_1)\in\mathbb{R}_{>0}^2$, the function $f_{t_0,t_1}$ is strictly increasing.
\end{lem01}
\begin{proof}
Let $(t_0,t_1)\in \mathbb{R}_{>0}^2$ be given. For any $s\geq 0$ we have
\[
f_{t_0,t_1}(s)
=
\min\left\{ \frac{s^2}{t_0s + (1-s)t_1}, \frac{s}{t_1}\right\}.
\]
Thus, it suffices to show that the functions $s\mapsto \frac{s}{t_1}$ and $s\mapsto \frac{s^2}{t_0s + (1-s)t_1}$ are strictly increasing. This is obvious for $s\mapsto \frac{s}{t_1}$. When $t_0\leq t_1$, the function $s\mapsto t_0s + (1-s)t_1$ is non-increasing, so $s\mapsto \frac{s^2}{t_0s + (1-s)t_1}$ is strictly increasing. If $t_1<t_0$, then the derivative of
\[
s
\mapsto
\frac{s^2}{(t_0-t_1)s+t_1} 
\]
is positive for $s>0$ and it is $0$ for $s=0$. Therefore, the function $s
\mapsto \frac{s^2}{(t_0-t_1)s+t_1}$ is strictly increasing.
\end{proof}
For each $n\in\mathbb{N}$, $\mathbf{t}=(t_0,t_1)\in \mathbb{R}_{>0}^{2}$, and $B>1$, define the map $g_n^{\mathbf{t}}(\,\cdot\,;B):\mathbb{R}_{\geq 0}\to\mathbb{R}_{>0}$ by
\[
g_n^{\mathbf{t}}(\rho;B) \colon= \sum_{d=2}^n \frac{1}{d^{\rho}(d-1)^{\rho}B^{f_{t_0,t_1}(\rho)}}
\;
\text{ for all }\rho\in \mathbb{R}_{>0}
\]
and $g^{\mathbf{t}}(\,\cdot\,;B):\mathbb{R}_{\geq 0}\to \mathbb{R}_{>0}\cup\{\infty\}$ by
\[
g(\rho;B) \colon= \sum_{d=2}^{\infty} \frac{1}{d^{\rho}(d-1)^{\rho}B^{f_{t_0,t_1}(\rho)}}
\;
\text{ for all }\rho\in \mathbb{R}_{>0}.
\]
Observe that for all $B>1$ we have
\[
g^{\mathbf{t}}(1;B) 
= 
\sum_{d=2}^{\infty} \frac{1}{d(d-1)B^{1/t_1}}
<
\sum_{d=2}^{\infty} \frac{1}{d(d-1)}=1
\]
and
\[
g^{\mathbf{t}}\left(\frac{1}{2};B\right)
=
\sum_{d=2}^{\infty} \frac{1}{d^{1/2}(d-1)^{1/2}B^{f_{t_0,t_1}(1/2)}}
=
\frac{1}{B^{f_{t_0,t_1}(1/2)}}
\sum_{d=2}^{\infty} \frac{1}{d^{1/2}(d-1)^{1/2}}
=
\infty,
\]
hence $\frac{1}{2} < s_0(B) <1$. Also, every $m,n\in\mathbb{N}$ with $m<n$ and every $\rho>0$ satisfy
\[
g_m^{\mathbf{t}}(\rho;B)
< 
g_n^{\mathbf{t}}(\rho;B)
<
g^{\mathbf{t}}(\rho;B).
\]
As a consequence, the sequence $(s_n(B))_{n\geq 1}$ is strictly increasing and each of its terms is bounded above by $s_0(B)$.
\begin{lem01}\label{LEM:sns0}
We have 
\[
\lim_{n\to\infty} s_n(B) = s_0(B).
\]
\end{lem01}
\begin{proof}
The discussion preceding the lemma implies
\[
\lim_{n\to\infty} s_n(B)\leq s_0(B).
\]
Pick any positive number $t<s_0(B)$. Then, $g(t;B)>1$ and every large $n\in \mathbb{N}$ satisfies $g_n(t;B)>1$, so $t< s_n(B)$ and 
\[
s_0(B)\leq \lim_{n\to\infty} s_n(B).
\]
\end{proof}
\begin{lem01}\label{LEM:s0:cont}
The function $B\mapsto s_0(B)$ is continuous.
\end{lem01}
\begin{proof}
Fix $B>1$. Take $n\in\mathbb{N}$. Let $\veps>0$ be arbitrary and
\[
0
<
\delta
=
\frac{1}{2}
\min
\left\{ 
B - B^{\frac{f_{t_0,t_1}(s_n(B))}{f_{t_0,t_1}(s_n(B)+\veps )}},
B^{\frac{f_{t_0,t_1}(s_n(B))}{f_{t_0,t_1}(s_n(B) - \veps  )}} - B
\right\}.
\]
Then, we have
\begin{align*}
g_n(s_n(B)+\veps;B-\delta) 
&=
\sum_{d=2}^n \frac{1}{d^{s_n(B)+\veps}(d-1)^{s_n(B)+\veps}(B - \delta)^{f_{t_0,t_1}(s_n(B)+\veps)}} \\
&\leq
\frac{B^{f_{t_0,t_1}(s_n(B))} }{(B - \delta)^{f_{t_0,t_1}(s_n(B)+\veps)}}
\sum_{d=2}^n \frac{1}{d^{s_n(B)}(d-1)^{s_n(B)}B^{f_{t_0,t_1}(s_n(B))}} \\
&=
\frac{B^{f_{t_0,t_1}(s_n(B))} }{(B-\delta)^{f_{t_0,t_1}(s_n(B)+\veps)}} < 1.
\end{align*}
Since $s_n$ is non-increasing on $B$, we conclude
\[
s_n(B)
\leq 
s_n(B-\delta)
<
s_n(B)+\veps.
\]
Similarly, we can show $g_n(s_n(B)-\veps;B+\delta)>1$ and, hence, 
\[
s_n(B+\delta)
\leq 
s_n(B)
<
s_n(B+\delta) + \veps.
\]
As a consequence, $s_n$ is continuous.

Since $s_n(B)\to s_0(B)$ as $n\to\infty$ and since $f_{t_0,t_1}$ is continuous and strictly increasing, we may pick $\delta>0$ such that every large $m\in\mathbb{N}$ verifies
\[
0
<
\delta
=
\frac{1}{2}
\min
\left\{ 
B - B^{\frac{f_{t_0,t_1}(s_m(B))}{f_{t_0,t_1}(s_m(B) + \veps )}},
B^{ \frac{f_{t_0,t_1}(s_m(B))}{ f_{t_0,t_1}(s_m(B) - \veps )}} - B
\right\}.
\]
For such $m\in\mathbb{N}$, every $B'\in \mathbb{R}_{>1}$ with the property $|B-B'|<\delta$ satisfies $| s_m(B) - s_m(B') |<\veps$. Letting $m\to \infty$, we conclude $| s_0(B) - s_0(B')|\leq \veps$. Therefore, $s_0$ is continuous.
\end{proof}
Our argument actually proves the next result.
\begin{teo01}
Let $f:\mathbb{R}_{>0}\to\mathbb{R}_{>0}$ be a strictly increasing continuous function. For each $n\in\mathbb{N}$, let $s=s_n(B,f)$ be the unique solution $s$ of 
\[
\sum_{d=2}^n \frac{1}{d^s(d-1)^sB^{f(s)}}=1.
\]
Call $s_0(B,f)$ the unique solution $s$ of
\[
\sum_{d=2}^{\infty} \frac{1}{d^s(d-1)^sB^{f(s)}}=1.
\]
Then, $s_n(B,f)\to s_0(B,f)$ as $n\to\infty$.
\end{teo01}

\subsection{Hausdorff dimension estimates}
We split the proof of Theorem \ref{LEM:TEO:Aux} into two parts. First, we use a particular family of coverings of $\clE_{\mathbf{t}}(\Psi)$ to obtain an upper bound for $\dimh \clE_{\mathbf{t}}(\Psi)$. The lower bound is proved considering two cases. In one of them, we use Lemma \ref{LEM:WWB-LUROTH}. In the other case, we apply the Mass Distribution Principle (Lemma \ref{LE:MDP}). 

The hypothesis $B>1$ implies $\liminf_n \Psi(n)= \infty$. Then, without loss of generality, we may assume that
\[
\Psi(n)\geq 2^{mt} 
\quad\text{ for all } n\in\mathbb{N}.
\]
\subsubsection{Upper bound}
\begin{proof}[Proof of Lemma \ref{LEM:TEO:Aux}. Upper bound.]
For each real number $A$ satisfying $1<A<B$, define the sets
\begin{align*}
\clE_{\mathbf{t}}'(A)
&\colon= 
\left\{ x\in (0,1]: d_n^{t_0}(x) \leq A^n \text{ and } d_{n+1}^{t_1}(x) \geq \frac{B^n}{d_n^{t_0}(x)} \text{ for infinitely many } n\in\mathbb{N} \right\}, \\
\clE_{\mathbf{t}}''(A)
&\colon= 
\left\{ x\in (0,1]: d_n^{t_0}(x) \geq A^n \text{ for infinitely many } n\in\mathbb{N} \right\}.
\end{align*}
Then, $\clE_{\mathbf{t}}(B) \subseteq \clE_{\mathbf{t}}'(A)\cup \clE_{\mathbf{t}}''(A)$, so
\[
\dimh\clE_{\mathbf{t}}(B) \leq \max\left\{ \dimh \clE_{\mathbf{t}}'(A), \dimh \clE_{\mathbf{t}}''(A)\right\}.
\]
Lemma \ref{LEM:WWB-LUROTH} implies
\[
\dimh \clE_{\mathbf{t}}''(A)=s\left(A^{1/t_0}\right).
\]
We take advantage of $\clE_{\mathbf{t}}'(A)$ being a limsup-set to give an upper estimate for its dimension. For each $\mathbf{d}=(d_1,\ldots,d_n)\in \scD^n$, define
\[
J_n(\mathbf{d})
\colon=
\bigcup
\left\{ I_{n+1} (\mathbf{d} d_{n+1}) : d_{n+1} \geq \max\left\{2, \left(\frac{B^{n}}{d_{n}^{t_0}} \right)^{1/t_1} \right\}\right\}
\]
Observe that $d_n^{t_0}\leq A^n$ and $d_{n+1}^{t_1} \geq B^n/d_n^{t_0}$ imply 
\[
d_{n+1}
\geq 
\left( \frac{B}{A}\right)^{n/t_1},
\]
and $(B/A)^{n/t_1}> 2$ for every large $n$ (depending on $B$, $A$, and $t_1$). For such $n\in\mathbb{N}$, we have
\begin{align*}
|J_n(\mathbf{d})|
&=
\left(\prod_{j=1}^n \frac{1}{d_j(d_j-1)}\right) \frac{1}{\left\lceil \frac{B^n}{d_{n}^{t_0/t_1}}\right\rceil -1 } \\
&\asymp
\left(\prod_{j=1}^n \frac{1}{d_j(d_j-1)}\right) \frac{d_{n}^{t_0/t_1}}{B^{n/t_1}} \\
&\asymp
\frac{1}{d_n^{2 - \frac{t_0}{t_1}}B^{\frac{n}{t_1}} }
\prod_{j=1}^{n-1} \frac{1}{d_j(d_j-1)}
=
\frac{|I_{n-1}(d_1,\ldots,d_{n-1})|}{d_n^{2 - \frac{t_0}{t_1}}B^{\frac{n}{t_1}}}.
\end{align*}
Take any $s>0$ and let $\clH^s$ be the $s$-Hausdorff measure on $\mathbb{R}$. From the inequality
\[
\clE'_{\mathbf{t}}(A)
=
\bigcap_{N=1}^{\infty} \bigcup_{n=N}^{\infty}\bigcup_{\mathbf{d}\in\scD^{n-1}} \bigcup_{d_n=2}^{\lfloor A^{n/t_0}\rfloor} J_n(\mathbf{d} \, d_n)
\]
we conclude that
\begin{align*}
\clH^s(\clE_{\mathbf{t}}'(A))
&\leq
\liminf_{N\to\infty} \sum_{n\geq N} \sum_{\mathbf{d}\in\scD^{n-1}} \sum_{d_n=2}^{\lfloor A^{n/t_0} \rfloor} |J_{n}(\mathbf{d}\, d_n)|^s \\
&\asymp
\liminf_{N\to\infty} \sum_{n\geq N} \sum_{\mathbf{d}\in\scD^{n-1}} \sum_{d_n=2}^{\lfloor A^{n/t_0} \rfloor}  |I_{n-1}(\mathbf{d})|^s d_n^{-s\left(2 - \frac{t_0}{t_1}\right)} B^{ -\frac{sn}{t_1} } \\
&=
\liminf_{N\to\infty} \sum_{n\geq N} \sum_{\mathbf{d}\in\scD^{n-1}} \frac{|I_{n-1}(\mathbf{d})|^s}{ B^{ \frac{sn}{t_1} }} \sum_{d_n=2}^{\lfloor A^{n/t_0} \rfloor}    d_n^{-s\left(2 - \frac{t_0}{t_1}\right)} \\
&=
\liminf_{N\to\infty} \sum_{n\geq N} 
\left(\sum_{ d\in\scD } \left( B^{ \frac{1}{t_1} } d(d-1)\right)^{-s}  \right)^n 
\sum_{d_n=2}^{\lfloor A^{n/t_0} \rfloor}    d_n^{-s\left(2 - \frac{t_0}{t_1}\right)} \\
&\ll
\liminf_{N\to\infty} \sum_{n\geq N} 
\left(\sum_{ d\in\scD } \left( B^{ \frac{1}{t_1} } d(d-1)\right)^{-s}  \right)^n 
\max
\left\{ 1,A^{\frac{n}{t_0} \left( 1 - s \left(2-\frac{t_0}{t_1} \right)\right)}\right\} \\
&=
\liminf_{N\to\infty} \sum_{n\geq N} 
\left(\sum_{ d\in\scD } \left( B^{ \frac{1}{t_1} } d(d-1)\right)^{-s}\max\left\{ 1,A^{\frac{1}{t_0} \left( 1 - s \left(2-\frac{t_0}{t_1} \right)\right)}\right\}  \right)^n .
\end{align*}
Therefore, if $S(A)$ is the solution of 
\[
\sum_{ d=2 }^{\infty} \left( B^{ \frac{1}{t_1} } d(d-1)\right)^{-s}\max\left\{ 1,A^{\frac{1}{t_0} \left( 1 - s \left(2-\frac{t_0}{t_1} \right)\right)}\right\}=1,
\]
we conclude that $\dimh \clE_{\mathbf{t}}(B)\leq S(A)$. Let $A>1$ be such that $S(A)=s(A)$, which occurs precisely when
\[
\frac{1}{A^{s/t_0}}
=
\frac{\max\{1,A^{(1-s(2-t_0/t_1))/t_0} \}}{B^{s/t_1}},
\]
or equivalently
\[
-f_{t_0}(s)\log A
=
\max\left\{0, \frac{1-2s}{t_0} + \frac{s}{t_1}\right\}\log A - f_{t_1}(s)\log B,
\]
which is
\[
f_{t_0,t_1}(s)\log B =  f_{t_0}(s) \log A.
\]
Using the exact same argument as in \cite[Lemma 5.1]{BakHusKleWan2022}, we conclude $1<A<B$ and, therefore, $\dimh \clE_{\mathbf{t}}(B)\leq s(A)=s_0(B)$.
\end{proof}

\subsubsection{Lower bound}
We consider two cases:
\[
\frac{s_0(B)}{t_1} - \frac{2s_0(B)-1}{t_0}\leq 0
\quad\text{ and }\quad
\frac{s_0(B)}{t_1} - \frac{2s_0(B)-1}{t_0}> 0.
\]
\subsubsection{Lower bound: first case}
Assume that $\frac{s_0(B)}{t_1} - \frac{2s_0(B)-1}{t_0}\leq 0$.
\begin{lem01}\label{LE:LEM:TEO:Aux:01}
If $\frac{s_0(B)}{t_1} - \frac{2s_0(B)-1}{t_0}\leq 0$, then $s = s_0(B)$ is the unique solution of 
\[
\sum_{d=2}^{\infty} \left( \frac{1}{d(d-1) B^{1/t_1}}\right)^s = 1.
\]
\end{lem01}
\begin{proof}
Since $f_{t_0,t_1}(s)\leq \frac{s}{t_1}$ for all $s>0$, we have 
\[
\sum_{d=2}^{\infty} \left( \frac{1}{d(d-1) B^{1/t_1}}\right)^s
\leq  
\sum_{d=2}^{\infty} \frac{1}{d^s(d-1)^sB^{f_{t_{0},t_1}(s) }}.
\]
Therefore, for all $\veps>0$, the inequality
\[
\sum_{d=2}^{\infty} \left( \frac{1}{d(d-1) B^{1/t_1}}\right)^{s_0(B)+\veps}
<1,
\]
holds and we conclude $s(B^{1/t_1})\leq s_0(B)$. 

We consider two further sub-cases. First, if $\frac{s_0(B)}{t_1} - \frac{2s_0(B) - 1}{t_0}< 0$, then $f_{t_0,t_1}(s)=\frac{s}{t_1}$ for every $s$ sufficiently close to $s_0$ and for any sufficiently small $\veps>0$ we have 
\[
1
<
\sum_{d=2}^{\infty} \left( \frac{1}{d(d-1) B^{1/t_1}}\right)^{s_{0}(B) - \veps }
=
\sum_{d=2}^{\infty} \frac{1}{d^s(d-1)^sB^{f_{t_{0},t_1}(s_{0}(B) - \veps ) }}.
\]
This implies $s_0(B)\leq s(B^{1/t_1}) + \veps$ and, letting $\veps\to 0$, we get $s_0(B)\leq s(B^{1/t_1})$. Second, assume that $\frac{s_0(B)}{t_1} - \frac{2s_0(B)-1}{t_0} =0$. Let $\delta'>0$ be arbitrary. Since $f_{t_0,t_1}$ is continuous, for any $\veps>0$, there is $\delta''>0$ such that, for all $s\in\mathbb{R}$,
\[
s_0(B)-\delta'' < s < s_0(B)
\quad\text{ implies }\quad
f_{t_0,t_1}(s)> \frac{s_0(B)}{t_1} - \veps.
\]
As a consequence, for $\delta=\min\{\delta',\delta''\}$ and any $s$ such that $s_0(B)-\delta < s < s_0(B)$, we have
\[
1 
\leq 
\sum_{d=2}^{\infty} \frac{1}{d^s(d-1)^sB^{f_{t_0,t_1}(s)}}
\leq 
\sum_{d=2}^{\infty} \frac{1}{d^s(d-1)^sB^{\frac{s_0(B)}{t_1} - \veps }}.
\]
Hence, $s(B^{1/t_1}) -\delta\geq s_0(B)$. Letting $\delta'\to 0$, we have $\delta\to 0$ and, thus, $s_0(B) \leq s(B^{1/t_1})$.
\end{proof}
\begin{lem01}\label{LE:LEM:TEO:Aux:02}
If $\frac{s_0(B)}{t_1} - \frac{2s_0(B)-1}{t_0}\leq 0$, then $\dimh \clE_{\mathbf{t}}(B) \geq s_0(B)$.
\end{lem01}
\begin{proof}
The definition of $\clE_{\mathbf{t}}(B)$ yields
\[
\left\{ x\in (0,1]: d_{n+1}^{t_1}(x)\geq B^n \text{ for infinitely many } n\right\}
\subseteq
\clE_{\mathbf{t}}(B).
\]
Then, by Lemmas \ref{LEM:WWB-LUROTH} and \ref{LE:LEM:TEO:Aux:01}, we conclude $\dimh \clE_{\mathbf{t}}(B) \geq s(B^{1/t_1}) = s_0(B)$.
\end{proof}

\subsubsection{Lower bound: second case}
Assume that
\begin{equation}\label{EQ:TEO:MDP:B}
\frac{s_0(B)}{t_1} - \frac{2s_0(B)-1}{t_0}> 0.
\end{equation}
In what follows, for all $x\in \mathbb{R}$ and all $r>0$, we write $B(x;r)\colon=(x-r,x+r)$.
\begin{lem01}[Mass Distribution Principle]\label{LE:MDP}
Let $F\subseteq \mathbb{R}$ be a non-empty set and let $\mu$ be a finite measure satisfying $\mu(F)>0$. If there are constants $c>0$, $r_0>0$, and $s\geq 0$ such that
\[
\mu\left( B(x;r)\right) \leq c r^s
\;\text{ for all } \;  x\in F \;\text{ and }\; r\in (0,r_0),
\]
then $\dimh F \geq s$.
\end{lem01}
\begin{proof}
see \cite[Proposition 2.1]{Fal1997}.
\end{proof}
\paragraph{Construction of the Cantor set.} 
By Lemma \ref{LEM:sns0} and \eqref{EQ:TEO:MDP:B}, we may take an $M\in\mathbb{N}_{\geq 3}$ such that $s\colon=s_M(B)$ is so close to $s_0(B)$ that
\[
\frac{1}{2}< s < 1
\quad\text{ and }\quad
f_{t_0,t_1}(s) = \frac{sf_{t_0}(s)}{t_1\left( f_{t_0}(s) + \frac{s}{t_1} - \frac{2s-1}{t_0}\right)}.
\]
Let $A$ be such that 
\[
f_{t_0,t_1}(s)\log A =  f_{t_0}(s) \log B;
\]
as above, $1<A<B$. Let $(\ell_{k})_{k\geq 1}$ be a sequence in $\mathbb{N}$ such that $
\ell_k \gg e^{\ell_1 + \ldots + \ell_{k-1}}$ for all $k\in\mathbb{N}$. Write
\[
\alpha_0\colon= A^{1/t_0}
\quad\text{ and }\quad
\alpha_1\colon= \left( \frac{B}{A}\right)^{1/t_1}
\]
and let $N\in\mathbb{N}$ be such that $2<\alpha_0^{N}$. Define the sequence $(n_j)_{j\geq 1}$ by
\[
n_1=\ell_1 N +1 
\quad\text{ and }\quad
n_{k+1} - n_k = \ell_{k+1} N + 2
\;\text{ for all } k\in \mathbb{N}.
\]
We can take the sequence $(\ell_k)_{k\geq 1}$ so sparse that $(n_k)_{k\geq 1}$ satisfies
\begin{equation}\label{EQ:LWEtB:00}
\frac{(2k-1)\log \alpha_0}{\log \alpha_1} < n_1 + \ldots + n_k
\quad\text{ for all } k\in\mathbb{N}_{\geq 2}.
\end{equation}
Let $E$ be the set formed by the real numbers $x=\langle d_1,d_2,d_3, \ldots\rangle$ satisfying the following conditions:
\begin{enumerate}[i.]
\item For every $k\in\mathbb{N}$, we have
\[
\alpha_0^{n_k} A^{n_k/t_0} \leq d_{n_k} \leq 2\alpha_0^{n_k}
\;\text{ and }\;
\alpha_1^{n_k}  \leq d_{n_k+1} \leq 2 \alpha_1^{n_k} .
\]
\item If $n\in \mathbb{N}\setminus\{n_k:k\in\mathbb{N}\}$, then $2\leq d_n\leq M$.
\end{enumerate}

Let us exhibit the Cantor structure of $E$. Define $D\colon= \Lambda^{-1}[E]$ (see Section \ref{SEC:LUROTH}) and for all $n\in\mathbb{N}$ define
\[
D_n
\colon= 
\left\{(d_1,\ldots, d_n)\in\scD^n: (d_j)_{j\geq 1} \in D\right\}.
\]
For each $\mathbf{d}\in D_n$ define the compact interval
\[
J_{n}(\mathbf{d})
\colon=
\bigcup_{\substack{d_{n+1}\in\scD\\ \mathbf{d} d\in D_{n+1}}} \overline{I}_{n+1}(\mathbf{d} d_{n+1}).
\]
We refer to the sets of the form $J_n(\mathbf{d})$ as \textbf{fundamental intervals} of order $n$. Clearly,
\[
E
=
\bigcap_{n\in\mathbb{N}} \bigcup_{\mathbf{d}\in D_n} J_n(\mathbf{d}).
\]
{\bf A probability measure.} Observe that for every $\mathbf{d}=(d_j)_{j\geq 1}\in D$ there is a finite collection of words $\bfw_{1}^{j}$, $\ldots$, $\bfw_{\ell_j}^{j}$ in $\{2,\ldots,M\}^N$ for $j\in\mathbb{N}$ such that, writing $\bfW_j\colon= \bfw_1^{j}\cdots \bfw_{\ell_j}^{j}$,
\begin{align*}
\mathbf{d}
&=
\bfw_{1}^{1}\bfw_{2}^{1}\ldots \bfw_{\ell_1}^{1} 
d_{n_{1}}d_{n_{1}+1} 
\bfw_{1}^{2}\bfw_{2}^{2}\ldots \bfw_{\ell_2}^{2} 
d_{n_{2}}d_{n_{2}+1} 
\cdots 
\bfw_{1}^{k}\bfw_{2}^{k}\ldots \bfw_{\ell_k}^{k} 
d_{n_{k}}d_{n_{k}+1} 
\cdots \\
&= 
\bfW_1
d_{n_{1}}d_{n_{1}+1} 
\bfW_2
d_{n_{2}}d_{n_{2}+1} 
\ldots 
\bfW_k
d_{n_{k}}d_{n_{k}+1} 
\ldots .
\end{align*}
Take $n\in\mathbb{N}$ and $\mathbf{d}\in D_n$. First, assume that $n\leq n_1+1$.
\begin{enumerate}[i.]
\item If $n=N\ell$ for some $\ell\in \{1,\ldots, \ell_1\}$, define
\[
\mu\left( J_{\ell N}(\mathbf{d})\right)
\colon=
\frac{|I_{N\ell}(\mathbf{d})|^s}{\alpha_0^{\ell Ns}}.
\]

\item If there is some $\ell\in \{1,\ldots, \ell_1\}$ such that $(\ell-1)N+1\leq n \leq \ell N-1$, put
\[
\mu\left( J_n(\mathbf{d}) \right)
\colon=
\sum_{\substack{\bfb\in \scD^{\ell N-n}\\\mathbf{d}\bfb\in D_{\ell N} }} \mu\left( J_{\ell N}(\mathbf{d} \bfb)\right).
\]

\item If $n=n_1$, then
\[
\mu\left( J_{n_1 }(\mathbf{d}) \right)
\colon=
\frac{\mu\left(J_{n_1-1}(d_1,\ldots, d_{n_{1}-1})\right)}{\lfloor 2\alpha_0^{n_1} \rfloor - \lceil \alpha_0^{n_1}\rceil}.
\]
\item If $n=n_1 + 1$, define
\[
\mu\left( J_{n_1 + 1}(\mathbf{d}) \right)
\colon=
\frac{\mu\left(J_{n_1-1}(d_1,\ldots, d_{n_{1}-1})\right)}{\left(\lfloor 2\alpha_1^{n_1} \rfloor - \lceil \alpha_1^{n-1}\rceil\right)\left(\lfloor 2\alpha_0^{n_1} \rfloor - \lceil \alpha_0^{n-1}\rceil\right)}
\]
\end{enumerate}
Assume that we have defined $\mu$ on the fundamental sets of level $1,\ldots, n_{k}+1$ for some $k\in\mathbb{N}$. Suppose that $n\in\{n_{k}+2, \ldots, n_{k+1}+1\}$.
\begin{enumerate}[i.]
\item If $n=n_k+1 + N\ell$ for some $\ell\in \{1,\ldots, \ell_{k+1}\}$, write
\[
\mu\left( J_{n_k+1 +N\ell}(\mathbf{d})\right)
\colon=
\mu\left( J_{n_k+1}(d_1,\ldots, d_{n_k+1}) )\right) 
\frac{\left| I_{N\ell}(\bfw_{1}^{k+1} \cdots \bfw_{\ell}^{k+1} )\right|^s}{\alpha_0^{sN\ell}}.
\]
\item When $n_k+1 + (\ell-1)N< n < n_kj+1+\ell N$ for some $\ell\in\{1,\ldots, \ell_{k+1}\}$, define
\[
\mu\left( J_{n}(\mathbf{d})\right)
\colon=
\sum_{ \bfb } \mu\left( J_{n_k+1 + \ell N}(\mathbf{d} \bfb)\right),
\]
where the sum runs along those words $\bfb$ such that $\mathbf{d}\bfb\in D_{n_k+1 + \ell N}$.
\item When $n=n_{k +1}$, put
\[
\mu\left( J_{n_{k +1} }(\mathbf{d})\right)
\colon=
\frac{\mu\left( J_{n_{k +1} - 1}(\mathbf{d})\right)}{ \lfloor 2\alpha_0^{ n_{k+1} } \rfloor - \lceil \alpha_0^{n_{k+1}}\rceil }.
\]
\item Whenever $n=n_{k+1} + 1$, put
\[
\mu\left( J_{n_{k +1} + 1}(\mathbf{d})\right)
\colon=
\frac{\mu\left(J_{n_{k+1}-1}(d_1,\ldots, d_{n_{k+1}-1})\right)}{\left(\lfloor 2\alpha_1^{n_{k+1}} \rfloor - \lceil \alpha_1^{n_{k+1}}\rceil\right)\left(\lfloor 2\alpha_0^{n_{k+1}} \rfloor - \lceil \alpha_0^{n_{k+1}}\rceil\right)}
\]
\end{enumerate}
The procedure defines a probability measure on the fundamental sets of a given level. The choice of $A$ and $B$ and the definition of $\mu$ ensure the consistency conditions. Hence, by the Daniell-Kolmogorov Consistency Theorem \cite[Theorem 8.23]{Kal2021}, the function $\mu$ is indeed a probability measure on $E$.
\paragraph{Gap estimates.} For $n\in\mathbb{N}$ and $\mathbf{d}\in D_n$, let $G_n(\mathbf{d})$ be the distance between $J_n(\mathbf{d})$ and the fundamental interval of level $n$ closest to it; that is,
\[
G_{n}(\mathbf{d})
\colon=
\inf\left\{ d\left( J_n(\mathbf{d}), J_n(\bfe)\right): \bfe\in D_n, \, \bfe\neq \mathbf{d}\right\}.
\]
\begin{lem01}\label{Lem:Gap}
For any $n\in \mathbb{N}$ and any $\mathbf{d}\in D_n$, we have
\[
G_n(\mathbf{d})\geq \frac{1}{M}|I_n(\mathbf{d})|.
\]
\end{lem01}
\begin{proof}
The proof is by induction on $n$. Pick $d\in \{2,\ldots, M\}$. If $2\leq d\leq M-1$, then
\[
\inf J_1(d) - \sup J_1(d+1) 
= 
\inf J_1(d) - \sup I_1(d+1) 
=
\frac{\left| I_1(d)\right|}{M}.
\]
If $3\leq d\leq M$, then
\[
\sup I_1(d) - \inf J_1(d-1) =  \frac{1}{M}\left| I_1(d-1)\right| > \frac{1}{M}\left| I_1(d)\right|.
\]
This shows the result for $n=1$. Assume that the lemma holds for $n=\tilde{n}-1\in\mathbb{N}$. Suppose that either 
\begin{equation}\label{EQ:LEM:GAP:01}
1<\tilde{n}\leq N \ell_1 -1 
\;\text{ or }\;
n_k+1\leq \tilde{n}\leq n_k+1 + (\ell_k N -1) 
\;\text{ for some }k\in\mathbb{N}
\end{equation}
and consider $\mathbf{d}=(d_1,\ldots, d_{\tilde{n}})\in D_{\tilde{n}}$. When $2\leq d_{\tilde{n}} \leq M-1$, we have
\begin{align*}
\inf J_{\tilde{n}}(\mathbf{d}) - \sup J_{\tilde{n}}\left(d_1,\ldots, d_{\tilde{n}}+1\right)
&=
\inf J_{\tilde{n}}(\mathbf{d}) - \sup I_{\tilde{n}}\left(d_1,\ldots, d_{\tilde{n}}+1\right) \\
&=
\inf J_{\tilde{n}}(\mathbf{d}) - \inf I_{\tilde{n}}\left(\mathbf{d}\right) \\
&=
\frac{1}{M}\left| I_{\tilde{n}}\left(\mathbf{d}\right)\right|.
\end{align*}
If $3\leq d_{\tilde{n}} \leq M$, then 
\begin{align*}
\inf J_{\tilde{n}}(d_1,\ldots, d_{\tilde{n}} -1) - \sup J_{\tilde{n}}(\mathbf{d})
&=
\inf J_{\tilde{n}}(d_1,\ldots, d_{\tilde{n}} -1) - \sup I_{\tilde{n}}(\mathbf{d})\\
&=
\inf J_{\tilde{n}}(d_1,\ldots, d_{\tilde{n}} -1) - \inf I_{\tilde{n}}(d_1,\ldots, d_{\tilde{n}} -1)\\
&= 
\frac{1}{M} \left|I_{\tilde{n}}(d_1,\ldots, d_{\tilde{n}} -1)\right| \\
&>
\frac{1}{M} \left|I_{\tilde{n}}(\mathbf{d})\right|.
\end{align*}
We conclude that, provided $3\leq d_{\tilde{n}}\leq M-1$ or $\mathbf{d}=(2,2,\ldots,2)$, we have $G_{\tilde{n}}(\mathbf{d}) > M^{-1}|I_{\tilde{n}}(\mathbf{d})|$. Assume that $d_{\tilde{n}}=2$ and let $j\in\{1,\ldots, \tilde{n} - 1\}$ be the largest index such that $d_{j}\geq 3$. Then, the neighbor to the right of $J_{\tilde{n}}(\mathbf{d})$ is $J_{\tilde{n}}(d_1,\ldots, d_{j-1},M,\ldots,M)$ and, using the induction hypothesis on the second inequality, we have
\begin{align*}
\inf J_{\tilde{n}}(d_1,\ldots, d_{j-1},M,\ldots, M) - \sup J_{\tilde{n}}(\mathbf{d})
&> 
\inf J_{j}(d_1,\ldots, d_{j-1}) - \sup J_{j}(d_1,\ldots, d_j) \\
&>
\frac{1}{M}\left| I_{j}(d_1,\ldots, d_j)\right| \\
&>
\frac{1}{M}\left| I_{\tilde{n}}(\mathbf{d})\right|.
\end{align*}
A similar argument holds when $d_{\tilde{n}}=M$. This proves the result for $\tilde{n}$ assuming \eqref{EQ:LEM:GAP:01}.

Suppose that $\tilde{n}= n_k-1$ for some $k\in\mathbb{N}$. If $2\leq d_{\tilde{n}} \leq M-1$, then 
\begin{align*}
\inf J_{\tilde{n}}(\mathbf{d}) &- \sup J_{\tilde{n}}(d_1,\ldots, d_{\tilde{n}}+1 ) 
= \\ 
&=
\left(\inf J_{\tilde{n}}(\mathbf{d}) - \inf I_{\tilde{n}}(\mathbf{d})\right) + \left(\sup I_{\tilde{n}}(d_1,\ldots, d_{\tilde{n}}+1 ) - \sup J_{\tilde{n}}(d_1,\ldots, d_{\tilde{n}}+1 ) \right) \\
&=
\frac{1}{\lceil 2\alpha_0^{n_j}\rceil} \left| I_{\tilde{n}}(\mathbf{d})\right| 
+ 
\left( 1 - \frac{1}{\lfloor \alpha_0^{n_j}\rfloor -1 }\right)\left| I_{\tilde{n}}(d_1,\ldots, d_{\tilde{n}} +1)\right| \\
&>
\left( 1 - \frac{1}{\lfloor \alpha_0^{n_j}\rfloor -1 }\right)
\frac{d_{ \tilde{n}}-1 }{d_{ \tilde{n}} +1} 
\left| I_{\tilde{n}}(\mathbf{d})\right| 	\\
&=
\left( 1 - \frac{1}{\lfloor \alpha_0^{n_j}\rfloor -1 }\right)
\left( 1 - \frac{2 }{d_{ \tilde{n}} +1} \right)
\left| I_{\tilde{n}}(\mathbf{d})\right| 	\\
&>
\left( 1 - \frac{1}{\lfloor \alpha_0^{n_j}\rfloor -1 }\right)
\left( 1 - \frac{2 }{M} \right)
\left| I_{\tilde{n}}(\mathbf{d})\right|
>
\frac{1}{M}\left| I_{\tilde{n}}(\mathbf{d})\right|.
\end{align*}
The last inequality follows from $M>4$. A similar argument shows that, when $3\leq d_{\tilde{n}}\leq M$,
\[
\inf J_{\tilde{n}}(d_1,\ldots, d_{\tilde{n}}-1 ) - \sup J_{\tilde{n}}(\mathbf{d}) 
>
\frac{1}{M} \left| I_{\tilde{n}}(d_1,\ldots, d_{\tilde{n}}-1 )\right| \\
>
\frac{1}{M} \left| I_{\tilde{n}}(\mathbf{d})\right|.
\]
A slight modification of this argument yields the result for $\tilde{n}=n_k$.
\end{proof}
\paragraph{Length estimates.} 
For any $n\in \mathbb{N}$, let $c_n$ and $C_n$ be 
\[
c_n\colon=\min\{ d_n: (d_j)_{j\geq 1}\in D\}
\quad\text{ and }\quad
C_n\colon=\max\{ d_n: (d_j)_{j\geq 1}\in D\}.
\]
Hence, when $\mathbf{d}=(d_1,\ldots, d_n)\in D_n$, we have
\[
|J_n(\mathbf{d})|
= 
\left(\frac{1}{c_{n+1} - 1} - \frac{1}{C_{n+1}}\right)
\prod_{j=1}^{n}\frac{1}{d_j(d_j-1)}
=
\left(\frac{1}{c_{n+1} - 1} - \frac{1}{C_{n+1}}\right)
|I_{n}(\mathbf{d})|
\]
In particular, when $n\neq n_k-1$ and $n\neq n_k$ for all $k\in\mathbb{N}$, 
\[
|J_n(\mathbf{d})|
= 
\left(1 - \frac{1}{M}\right)
\prod_{j=1}^{n}\frac{1}{d_j(d_j-1)}
= 
\left(1 - \frac{1}{M}\right)
|I_{n}(\mathbf{d})|,
\]
which gives
\begin{equation}\label{Ec:LEST:01}
\frac{1}{2}
|I_{n}(\mathbf{d})|
\leq 
|J_n(\mathbf{d})|
\leq 
\left(1 - \frac{1}{M}\right)
|I_{n}(\mathbf{d})|.
\end{equation}
For $n=n_k-1$, we have
\[
|J_{n_k-1}(\mathbf{d})|
= 
\left(\frac{1}{\lceil \alpha_0^{n_k}\rceil -1} - \frac{1}{\lfloor 2 \alpha_0^{n_k}\rfloor}\right)
\prod_{j=1}^{n_k-1}\frac{1}{d_j(d_j-1)}
= 
\left(\frac{1}{\lceil \alpha_0^{n_k}\rceil -1} - \frac{1}{\lfloor 2 \alpha_0^{n_k}\rfloor}\right)
\left| I_{n_k-1}(\mathbf{d})\right|,
\]
so
\begin{equation}\label{Ec:LEST:02}
\frac{1}{2\alpha_0^{n_k}}
|I_{n_k-1}(\mathbf{d})|
\leq 
|J_{n_k-1}(\mathbf{d})|
\leq  
\frac{1}{\alpha_0^{n_k}}
|I_{n_k-1}(\mathbf{d})|.
\end{equation}
We can replace the constant $1$ in the upper bound with an arbitrary constant strictly larger than $\frac{1}{2}$, but we would have to consider larger values of $N$.

Similarly, for $n=n_k$,
\[
|J_{n_k}(\mathbf{d})|
=
\left(\frac{1}{\lceil \alpha_1^{n_k}\rceil -1} - \frac{1}{\lfloor 2 \alpha_1^{n_k}\rfloor}\right)
\prod_{j=1}^{n_k}\frac{1}{d_j(d_j-1)}
=
\left(\frac{1}{\lceil \alpha_1^{n_k}\rceil -1} - \frac{1}{\lfloor 2 \alpha_1^{n_k}\rfloor}\right)
\left| I_{n_k}(\mathbf{d})\right| 
\]
so
\[
\frac{1}{2\alpha_1^{n_k}}
|I_{n_k}(\mathbf{d})|
\leq 
|J_{n_k}(\mathbf{d})|
\leq  
\frac{1}{\alpha_1^{n_k}}
|I_{n_k}(\mathbf{d})|.
\]
Again, we can replace the constant $1$ in the upper bound with an arbitrary constant strictly larger than $\frac{1}{2}$ at the expense of a larger $N$.

\begin{lem01}\label{Lem510}
Let $k\in\mathbb{N}$ be arbitrary.
\begin{enumerate}[1.]
\item For any $\mathbf{d}=(d_1,\ldots, d_{n_k})\in D_{n_k}$, we have
\[
\frac{1}{8\alpha_0^{n_k}\alpha_1^{n_k}}
\left|J_{n_k-1}(d_1,\ldots, d_{n_k-1})\right|
<
\left| J_{n_k}(\mathbf{d})\right|.
\]
\item For any $\mathbf{d}=(d_1,\ldots, d_{n_k+1})\in D_{n_k+1}$, we have
\[
\frac{1}{2^6\alpha_0^{n_k}\alpha_1^{2n_k}}
\left|J_{n_k-1}(d_1,\ldots, d_{n_k-1})\right|
<
\left| J_{n_k +1 }(\mathbf{d})\right|.
\]
\end{enumerate}
\end{lem01}
\begin{proof}
\begin{enumerate}[1.]
\item By our previous discussion,
\begin{align*}
\left|J_{n_k}(\mathbf{d})\right|
&= 
\left(\frac{1}{\lceil \alpha_1^{n_k}\rceil -1} - \frac{1}{\lfloor 2 \alpha_1^{n_k}\rfloor}\right)
\left| I_{n_k}(\mathbf{d})\right| \\
&= 
\left(\frac{1}{\lceil \alpha_1^{n_k}\rceil -1} - \frac{1}{\lfloor 2 \alpha_1^{n_k}\rfloor}\right)
\frac{\left| I_{n_k-1}(d_1,\ldots, d_{n_k-1})\right|} {d_{n_k}(d_{n_k}-1)} \\
&> 
\frac{1}{8\alpha_1^{n_k} \alpha_0^{n_k}\alpha_0^{n_k}} \left| I_{n_k-1}(d_1,\ldots, d_{n_k-1})\right| \\
&> 
\frac{1}{8\alpha_1^{n_k} \alpha_0^{n_k}} \left|J_{n_k-1}(d_1,\ldots, d_{n_k-1})\right| && \text{(by \eqref{Ec:LEST:02})}.
\end{align*}
\item Similarly,
\begin{align*}
\left|J_{n_k +1 }(\mathbf{d})\right|
&=
\left(1 - \frac{1}{M}\right)
\left|I_{n_k +1}(\mathbf{d})\right| \\
&=
\left(1 - \frac{1}{M}\right)
\frac{\left|I_{n_k -1 }(d_1, \ldots, d_{n_k - 1} )\right|}{d_{n_k}(d_{n_k}-1)d_{n_k+1}(d_{n_k+1}-1) } \\ 
&>
\frac{1}{2^6} \frac{\left|I_{n_k - 1}(d_1, \ldots, d_{n_k - 1} )\right|}{\alpha_0^{2n_k}\alpha_1^{2n_k} } \\
&>
\frac{1}{2^6} \frac{\left|J_{n_k - 1}(d_1, \ldots, d_{n_k - 1} )\right|}{\alpha_0^{n_k}\alpha_1^{2n_k} } &&\text{(by \eqref{Ec:LEST:02})}.
\end{align*}
\end{enumerate}
\end{proof}

\paragraph{Measure of the fundamental intervals.} We now compute upper estimates of $\mu(J_n(\mathbf{d}))$ for all $n\in\mathbb{N}$ and $\mathbf{d}\in D_n$.
\begin{lem01}\label{Lem511}
The following statements hold:
\[
\frac{1}{\alpha_0} = \left( \frac{1}{\alpha_0\alpha_1}\right)^s
\quad\text{ and }\quad
\frac{1}{\alpha_0\alpha_1} \leq \left( \frac{1}{\alpha_0\alpha_1^2}\right)^s.
\]
\end{lem01}
\begin{proof}
See \cite[Lemma 6.6]{BakHusKleWan2022}.
\end{proof}

For each $k\in\mathbb{N}$, define the number $\gamma_{k}>0$ by
\[
\gamma_k^{-1} \colon=
\left( 1  - \frac{2}{\alpha_0^{n_k}}\right)
\left( 1  - \frac{2}{\alpha_1^{n_k}}\right).
\]
We will use the following obvious facts:
\begin{enumerate}[1.]
\item Since $\alpha_0>1$ and $\alpha_1>1$, the series $\sum_{k} \alpha_0^{-n_k}$ and  $\sum_{k} \alpha_0^{-n_k}$ are convergent and, thus, so is the product $C'\colon=\prod_{k\in\mathbb{N}} \gamma_k$ (see \cite[Theorem 7.32]{Str1981}).
\item For each $k\in\mathbb{N}$, 
\[
\frac{1}{\left(\lfloor 2\alpha_1^{n_{k-1}} \rfloor - \lceil \alpha_1^{n_{k-1}}\rceil\right)\left(\lfloor 2\alpha_0^{n_{k-1}} \rfloor - \lceil \alpha_0^{n_{k-1}}\rceil\right)}
< 
 \frac{\gamma_{k}}{\alpha_0^{n_k}\alpha_0^{n_k}}.
\]
\end{enumerate}
\begin{lem01}\label{PROPO:LWEtB:01}
For every $k\in\mathbb{N}_{\geq 2}$ and every $\mathbf{d}\in D_{n_k-1}$, we have 
\[
\mu\left( J_{n_k-1}(\mathbf{d})\right)
<
\gamma_{k-1} \left( \frac{1}{\alpha_0^{N\ell_k  }\alpha_1^{n_{k-1}} } \right)^s
\left( \frac{\left|I_{N\ell_k}(\bfW_k)\right| }{ \alpha_{0}^{n_{k-1}} \alpha_{1}^{n_{k-1}} }\right)^s
\mu\left( J_{n_{k-1}-1}(d_1,\ldots, d_{n_{k-1}-1})\right).
\]
\end{lem01}
\begin{proof}
Take $k$ and $\mathbf{d}$ as in the statement.  The lemma follows from the definition of $\mu$ and $n_k-1 = n_{k-1} + N\ell_k +1$:
\begin{align*}
\mu\left( J_{n_k-1}(\mathbf{d})\right)
&=
\frac{\left|I_{N\ell_k}\left( \bfW_{k}\right)\right|^{s}  }{ \alpha_0^{N\ell_ks } }
\mu\left( J_{n_{k-1}+1}(d_1,\ldots, d_{n_{k-1}+1})\right) \\
&=
\frac{\left|I_{N\ell_k}\left( \bfW_{k}\right)\right|^{s}  }{ \alpha_0^{N\ell_ks } }
\frac{\mu\left( J_{n_{k-1}-1}(d_1,\ldots, d_{n_{k-1}-1})\right)}{\left(\lfloor 2\alpha_1^{n_{k-1}} \rfloor - \lceil \alpha_1^{n_{k-1}}\rceil\right)\left(\lfloor 2\alpha_0^{n_{k-1}} \rfloor - \lceil \alpha_0^{n_{k-1}}\rceil\right)} \\
&<
\frac{\gamma_{k-1}}{\alpha_0^{N\ell_ks }}
\frac{\left|I_{N\ell_k}\left( \bfW_{k}\right)\right|^{s}  }{ \alpha_0^{n_{k-1}}\alpha_1^{n_{k-1}} }
\mu\left( J_{n_{k-1}-1}(d_1,\ldots, d_{n_{k-1}-1})\right) \\
&<
\frac{\gamma_{k-1}}{\alpha_0^{N\ell_ks }}
\frac{\left|I_{N\ell_k}\left( \bfW_{k}\right)\right|^{s}  }{ \alpha_0^{sn_{k-1}}\alpha_1^{2sn_{k-1}} }
\mu\left( J_{n_{k-1}-1}(d_1,\ldots, d_{n_{k-1}-1})\right).
\end{align*}
\end{proof}
\begin{lem01}\label{LE:TEO04:LB:BOLAS}
There exists a constant $C=C(B,M,N,s,\mathbf{t})>0$ such that 
\[
\mu\left( J_n(\mathbf{d}) \right) \leq C |J_n(\mathbf{d})|^{s}
\quad
\text{ for all } n\in\mathbb{N} \text{ and all } \mathbf{d}\in D_n.
\]
\end{lem01}
\begin{proof}
Pick $n\in\mathbb{N}$. First, we further assume that $1\leq n\leq n_1+1$.
\begin{enumerate}[i.]
\item Suppose that $n=\ell N$ for some $\ell\in \{1,\ldots, \ell_1\}$. When $1\leq \ell\leq \ell_1-1$, using \eqref{Ec:LEST:01} we get
\[
\mu\left( J_{\ell N}(\mathbf{d})\right)
=
\frac{1}{\alpha_0^{N\ell s}} \,
|I_{\ell N}(\mathbf{d})|^s 
=
\frac{2^s}{\alpha_0^{N\ell s}} \,
\frac{|I_{\ell N}(\mathbf{d})|^s}{2^s} 
<
\frac{2^s}{\alpha_0^{N\ell s}} \,
|J_{\ell N}(\mathbf{d})|^s 
<
|J_{\ell N}(\mathbf{d})|^s.
\]
When $\ell=\ell_1$, we have $n=\ell_1N = n_1 - 1$ and, by \eqref{Ec:LEST:02},
\[
\mu(J_{\ell_1N} (\mathbf{d}))
=
2^s\alpha_0^s \,
\frac{|I_{n_1-1}(\mathbf{d})|^s}{2^s\alpha_0^{sn_1}} 
<
2^s\alpha_0^s \,
|J_{\ell_1 N}(\mathbf{d})|^s
<
2 \alpha_0\,
|J_{\ell_1 N}(\mathbf{d})|^{s}.
\]
\item If $(\ell - 1)N + 1\leq n \leq \ell N-1$ for some $\ell\in\{1,\ldots, \ell_1\}$, then
\begin{align*}
\mu\left( J_n(\mathbf{d}) \right)
&=
\frac{1}{\alpha_0^{N\ell s}} \,
\sum_{\substack{\bfb\in \scD^{\ell N-n}\\\mathbf{d}\bfb\in D_{\ell N} }} \left| I_{\ell N}(\mathbf{d} \bfb)\right|^s \\
&=
\frac{|I_{n}(\mathbf{d})|^s}{\alpha_0^{ns}} \,
\sum_{\substack{\bfb\in \scD^{\ell N-n}\\\mathbf{d}\bfb\in D_{\ell N} }} \frac{\left| I_{\ell N - n}( \bfb)\right|^s}{\alpha_0^{(N\ell-n) s}} \\
&=
\frac{|I_{n}(\mathbf{d})|^s}{\alpha_0^{ns}} \,
\left( \sum_{d=2}^M \frac{1}{d^s(d-1)^s \alpha_0^{s}} \right)^{\ell N-n} \\
&=
\frac{|I_{n}(\mathbf{d})|^s}{\alpha_0^{ns}} \,
=
\frac{2^s}{\alpha_0^{ns}} \,
\frac{|I_{n}(\mathbf{d})|^s}{2^s}
<
\frac{2^s}{\alpha_0^{ns}} \,
|J_{n}(\mathbf{d})|^{s}
< |J_{n}(\mathbf{d})|^{s}.
\end{align*}
\item If $n=n_1$, then 
\begin{align*}
\mu\left( J_{n_1}(\mathbf{d})\right) 
&=
\frac{\mu\left( J_{n_1-1}(d_1,\ldots, d_{n_1-1} )\right) }{\lfloor 2\alpha_0^{n_1} \rfloor - \lceil \alpha_0^{n_1}\rceil } \\
&< 
\frac{2}{\alpha_0^{n_1}}  \mu\left( J_{n_1-1}(d_1,\ldots, d_{n_1-1} )\right) \\
&<
\frac{2}{\alpha_0^{n_1}}\,2\alpha_0|J_{n_1-1}(d_1,\ldots, d_{n_1-1})|^s \\
&=
2^2\, \alpha_0 8\frac{ |J_{n_1-1}( d_1,\ldots, d_{n_1-1} )|^s }{8 \left( \alpha_0^{n_1}\alpha_1^{n_1}\right)^s}
<
2^5 \alpha_0|J_{n_1}(\mathbf{d})|^s.
\end{align*}
We used Lemma \ref{Lem511} in the last inequality.
\item When $n=n_1+1$, we have 
\begin{align*}
\mu\left( J_{n_1+1}(\mathbf{d})\right) 
&< 
2^2\frac{\mu\left( J_{n_1+1}(d_1,\ldots, d_{n_k-1})\right) }{\alpha_0^{n_1}\alpha_1^{n_1}} \\
&\leq 
2^2
\frac{\mu\left( J_{n_1+1}(d_1,\ldots, d_{n_k-1})\right) }{\alpha_0^{n_1s}\alpha_1^{2n_1s}} \\
&<
2^3 \alpha_0\, 
\frac{\left| J_{n_1+1}(d_1,\ldots, d_{n_k-1})\right|^s }{\alpha_0^{n_1s}\alpha_1^{2n_1s}} 
<
2^9 \alpha_0\, 
\left| J_{n_1+1}(\mathbf{d} )\right|^s.
\end{align*}
\end{enumerate}
Assume now that $n_1+1<n$ and pick $k\in\mathbb{N}_{\geq 2}$ such that $n_k+1<n<n_{k+1}+1$.
\begin{enumerate}[i.]
\item If $n=n_{k}-1$, we apply Lemma \ref{PROPO:LWEtB:01} repeatedly to obtain
\begin{align*}
\mu\left( J_{n_k-1}(\mathbf{d})\right)
&<
\gamma_{k-1} \left( \frac{1}{\alpha_0^{N\ell_k  }\alpha_1^{n_{k-1}} } \right)^s
\left( \frac{\left|I_{N\ell_k}(\bfW_k)\right| }{ \alpha_{0}^{n_{k-1}} \alpha_{1}^{n_{k-1}} }\right)^s
\mu\left( J_{n_{k-1}-1}(d_1,\ldots, d_{n_{k-1}-1})\right) \\
&<
\gamma_1\cdots \gamma_{k-1} \left( \frac{1}{\alpha_0^{N\left(\ell_1 + \ldots + \ell_k\right)  }\alpha_1^{n_1 + \ldots + n_{k-1}} } \right)^s
\left( \frac{\left|I_{N\ell_1}(\bfW_1)\right| \cdots \left|I_{N\ell_k}(\bfW_k)\right| }{ \alpha_{0}^{n_{1}} \alpha_{1}^{n_{1}} \cdots  \alpha_{0}^{n_{k-1}} \alpha_{1}^{n_{k-1}} }\right)^s \\
&<
C' \left( \frac{ 1 }{\alpha_0^{N\left(\ell_1 + \ldots + \ell_k\right)  }\alpha_1^{n_1 + \ldots + n_{k-1}} } \right)^s
\left|I_{n_{k}-1 }(\mathbf{d})\right|^s \\
&\leq 
C' \left( \frac{ 2\alpha_0^{n_k}}{\alpha_0^{N\left(\ell_1 + \ldots + \ell_k\right)  }\alpha_1^{n_1 + \ldots + n_{k-1}} } \right)^s
\left|J_{n_{k}-1 }(\mathbf{d})\right|^s \\
&<
2
C' \left( \frac{1}{\alpha_0^{N\left(\ell_1 + \ldots + \ell_k\right) - n_k  }\alpha_1^{n_1 + \ldots + n_{k-1}} } \right)^s
\left|J_{n_{k}-1 }(\mathbf{d})\right|^s \\
&=
2
C' \left( \frac{\alpha_0^{2k-1}}{\alpha_1^{n_1 + \ldots + n_{k-1}} }  \right)^s
\left|J_{n_{k}-1 }(\mathbf{d})\right|^s  \\
&<
C'
\left|J_{n_{k}-1 }(\mathbf{d})\right|^s &\text{ (by \eqref{EQ:LWEtB:00}) }.
\end{align*}
\item Suppose we had $n=n_k$, then
\begin{align*}
\mu\left( J_{n_k}(\mathbf{d})\right)
&< 
\frac{2}{\alpha_{0}^{n_k}}
\mu\left( J_{n_k-1}(d_1,\ldots, d_{n_k - 1})\right) \\
&<
\frac{2}{\alpha_{0}^{n_k}}
C'
\left| J_{n_k-1}(d_1,\ldots, d_{n_k - 1})\right|^s \\
&<
2
C'
\frac{\left| J_{n_k-1}(d_1,\ldots, d_{n_k - 1})\right|^s}{\alpha_0^{2n_ks}\alpha_1^{2n_ks} } &&\text{(by Lemma \ref{Lem511})} \\
&< 
2^4
C'
\left| J_{n_k}(\mathbf{d} )\right|^s &&\text{(by Lemma \ref{Lem510}).}
\end{align*}
\item If $n=n_k+1$, then
\begin{align*}
\mu\left( J_{n_k+1}(\mathbf{d})\right)
&< 
\frac{2^2}{\alpha_{0}^{n_k}\alpha_{1}^{n_k}}
\mu\left( J_{n_k-1}(d_1,\ldots, d_{n_k - 1})\right) \\
&<
\frac{2^2}{\alpha_{0}^{n_k}\alpha_{1}^{n_k}}
C'
\left| J_{n_k-1}(d_1,\ldots, d_{n_k - 1})\right|^s \\
&<
2^2
C'
\left( \frac{\left| J_{n_k-1}(d_1,\ldots, d_{n_k - 1})\right|}{\alpha_0^{n_k}\alpha_1^{2n_k} } \right)^s \\
&<
2^2 \, 2^{6s} \,
C'
\left| J_{n_k+1}(\mathbf{d} )\right|^s &&\text{(by Lemma \ref{Lem511})} \\
&<
2^{8} \, C'
\left| J_{n_k+1}(\mathbf{d} )\right|^s.
\end{align*}
\item If $n=n_k + 1 +\ell N$ for some $1\leq \ell< \ell_{k+1}$, then
\begin{align*}
\mu\left( J_{n_k + 1 +\ell N}(\mathbf{d})\right)
&=
\frac{\left|I_{N\ell}(\bfw_1^{k+1}\cdots\bfw_{\ell}^{k+1})\right|^s}{\alpha_0^{N\ell s}}
\mu\left( J_{n_k + 1}(d_1,\ldots, d_{n_k+1} )\right) \\
&=
\frac{2^{8} \, C'}{\alpha_0^{N\ell s}}
\left|I_{N\ell}(\bfw_1^{k+1}\cdots\bfw_{\ell}^{k+1})\right|^s
\left| J_{n_k + 1}(d_1,\ldots, d_{n_k+1} )\right|^s \\
&<
2^{8} \, C'
\left| I_{n_k + 1}(d_1,\ldots, d_{n_k+1} )\right|^s
\frac{\left|I_{N\ell}(\bfw_1^{k+1}\cdots\bfw_{\ell}^{k+1})\right|^s}{\alpha_0^{N\ell s}} \\
&<
2^{8} \, C' \left| J_{n_k + 1 + N\ell}(\mathbf{d})\right|^s. 
\end{align*}
The last inequality is shown as in the case $k=1$.
\item Assume that $n_k + 2 + (\ell -1 ) N\leq n < n_k + 1 +\ell N$ with $1\leq \ell \leq \ell_k$, then
\begin{align*}
\mu(J_n(\mathbf{d}))
&< 
\mu\left( J_{n_k +1 + (\ell -1 ) N}(d_1,\ldots, d_{n_k +1 + (\ell -1 ) N})\right) \\
&<
2^{8} \, C'
\left| J_{n_k +1 + (\ell -1 ) N}(d_1,\ldots, d_{n_k +1 +(\ell -1 ) N})\right|^s.
\end{align*}
The discussion preceding Lemma \ref{Lem510} tells us that
\begin{align*}
\left| J_{n_k+1 + (\ell-1)N}(d_1,\ldots, d_{n_k+1 + (\ell-1)N}) \right| 
&\leq
\left( 1 - \frac{1}{M}\right) \left| I_{n_k+1 + (\ell-1)N}(d_1,\ldots, d_{n_k+1 + (\ell-1)N}) \right| \\
&=
\left( 1 - \frac{1}{M}\right) 
\left| I_n(\mathbf{d})\right|
\prod_{j= n_k+2 + (\ell-1)N}^n d_j(d_j-1) \\
&<
\left( 1 - \frac{1}{M}\right) 
M^{2N}\left| I_n(\mathbf{d})\right|  \\
&<
2\left( 1 - \frac{1}{M}\right) 
M^{2N}\left| J_n(\mathbf{d})\right|.
\end{align*}
As a consequence, we have
\[
\mu(J_n(\bfa))
< 
2^9C'\left( 1 - \frac{1}{M}\right) 
M^{2N}\left| J_n(\mathbf{d})\right|^s.
\]
\end{enumerate}
\end{proof}
{\bf Measure of balls.}
We now estimate the measure of balls with center on $E$ and small radius. Define
\[
r_0\colon= \min\left\{ G_1(d):2\leq d\leq M\right\}.
\]
\begin{lem01}
There is a constant $C''>0$ such that 
\[
\mu\left( B(x;r) \right) \leq C'' r
\quad
\text{ for all }x\in E \text{ and all } r\in (0,r_0).
\]
\end{lem01}
Take any $x=\langle d_1,d_2,\ldots\rangle$ and any $r\in (0,r_0)$. Pick $n\in\mathbb{N}$ such that
\[
G_{n+1}(d_1,\ldots, d_{n+1})
\leq 
r
<
G_n(d_1,\ldots,d_n).
\]
By definition of $G_n$, the ball $B(x;r)$ intersects exactly one fundamental interval of order $n$, namely $J_n(d_1,\ldots, d_n)$. 

Let us further assume that $n_k + 1 \leq n \leq n_{k+1}-1$ for some $k\in\mathbb{N}$. Taking $C$ as in Lemma \ref{LE:TEO04:LB:BOLAS},
\begin{align*}
\mu\left( B(x;r)\right)
&\leq 
\mu\left(J_n(d_1,\ldots, d_n) \right) \\
&\leq 
C
\left|J_{n}(d_1,\ldots, d_{n}) \right|^s \\
&=
C
\left|I_{n}(d_1,\ldots, d_{n}) \right|^s\left(\frac{M-1}{M}\right)^s  \\
&=
C
\left|I_{n+1}(d_1,\ldots, d_{n},d_{n+1}) \right|^s\left(\frac{M-1}{M}\right)^s d_{n+1}^s(d_{n+1}-1)^s  \\
&<
CM^{3}
G_{n+1}(d_1,\ldots, d_{n+1})^s &&\text{(by lemma \eqref{Lem:Gap})} \\
&<
CM^{3} r^s.
\end{align*}
Suppose now that $n=n_k-1$. We consider two cases: 
\begin{equation}\label{EQ:LWEtB:01}
r< \frac{\left| I_{n_k-1}(d_1,\ldots, d_{n_k-1})\right|}{\alpha_0^{2n_k}}
\quad\text{ and }\quad
r \geq \frac{\left| I_{n_k-1}(d_1,\ldots, d_{n_k-1})\right|}{\alpha_0^{2n_k}}.
\end{equation}
In the first case, $B(x;r)$ intersects at most three fundamentals intervals of level $n_k$. As a consequence, we have
\begin{align*}
\mu\left( B(x;r)\right)
&\leq 
3 \mu\left( J_{n_k}(d_1,\ldots, d_{n_k})\right)\\
&\leq
3 C \left| J_{n_k}(d_1,\ldots, d_{n_k})\right|^s 
\leq
3 C M  G_{n_k}(d_1,\ldots, d_{n_k})^s
\leq
3 C M r^s.
\end{align*}
Assume the second inequality in \eqref{EQ:LWEtB:01}. All the cylinders of level $n_k$ contained in $J_{n_j-1}(d_1,\ldots, d_{n_k-1})$ are of the form 
\[
I_{n_k}(d_1,\ldots, d_{n_k-1},a)
\quad\text{ with }\quad
a\in\{ \lceil \alpha_0^{n_k}\rceil, \ldots, \lfloor 2 \alpha_0^{n_k}\rfloor\},
\]
so
\[
\left| I_{n_k}(d_1,\ldots, d_{n_k-1},a)\right|
=
\frac{\left| I_{n_k-1}(d_1,\ldots, d_{n_k-1})\right|}{a(a-1)}
\geq 
\frac{\left| I_{n_k-1}(d_1,\ldots, d_{n_k-1})\right|}{\alpha_0^{2n_k}};
\]
If $T$ is the total amount of cylinders of level $n_k$ contained in $B(x;r)$, then
\[
T
\leq 
\frac{\alpha_0^{2n_k}}{ \left| I_{n_{k-1}}(d_1,\ldots, d_{n_{k-1}})\right|}
\, 2r
\]
and the total amount of cylinders of level $n_k$ intersecting $B(x;r)$ is at most
\[
\frac{2r\alpha_0^{2n_k}}{ \left| I_{n_{k-1}}(d_1,\ldots, d_{n_{k-1}})\right|}
+ 2
<
\frac{4r\alpha_0^{2n_k}}{ \left| I_{n_{k-1}}(d_1,\ldots, d_{n_{k-1}})\right|}.
\]
Since each cylinder of level $n_k$ contains at most one fundamental interval of level $n_k$, we have
\begin{align*}
\mu\left( B(x;r)\right)
&< 
\frac{4r\alpha_0^{2n_k}}{ \left| I_{n_{k-1}}(d_1,\ldots, d_{n_{k-1}})\right|} \mu\left( J_{n_k}(d_1,\ldots, d_{n_k-1},d_{n_k})\right) \\
&<
\frac{8r\alpha_0^{2n_k} }{ \left| I_{n_{k-1}}(d_1,\ldots, d_{n_{k-1}})\right|} \mu\left( J_{n_k-1}(d_1,\ldots, d_{n_k-1}) \right).
\end{align*}
The second inequality follows from the definition of $\mu\left( J_{n_k}(d_1,\ldots, d_{n_k-1},d_{n_k})\right)$. Take $C$ as in Lemma \ref{LE:TEO04:LB:BOLAS}. Then, since $\min\{a,b\} \leq a^{1-s} b^s$ for all positive $a,b$, 
\begin{align*}
\mu\left( B(x;r)\right)
&\leq
\min\left\{ 
\mu(J_{n_{k-1}}(d_1,\ldots, d_{n_{k-1}})), 
8r\alpha_0^{2n_k} \,
\frac{\mu\left( J_{n_k-1}(d_1,\ldots, d_{n_k-1}) \right)}{ \left| I_{n_{k-1}}(d_1,\ldots, d_{n_{k-1}})\right|} 
\right\} \\
&=
\mu(J_{n_{k-1}}(d_1,\ldots, d_{n_{k-1}}))
\min\left\{ 
1, 
\frac{ 8r\alpha_0^{2n_k} }{ \left| I_{n_{k-1}}(d_1,\ldots, d_{n_{k-1}})\right|} 
\right\} \\
&<
8^sr^s\alpha_0^{2sn_k} \,
\frac{\mu(J_{n_{k-1}}(d_1,\ldots, d_{n_{k-1}}))}{\left| I_{n_{k-1}}(d_1,\ldots, d_{n_{k-1}})\right|^s} \\
&<
8Cr^s.
\end{align*}
In the last inequality, we have used \eqref{Ec:LEST:02}. A similar argument holds for $n=n_k$, since we have distributed uniformly the mass $\mu(J_{n_{k-1}}(d_1,\ldots, d_{n_k-1} ))$ among the $\lfloor 2 \alpha_1^{n_k}\rfloor - \lceil \alpha_1^{n_k}\rceil$ fundamental intervals of level $n_k$ contained in $J_{n_{k-1}}(d_1,\ldots, d_{n_k-1})$ when defining $\mu$.
\begin{proof}[Proof of Theorem  \ref{TEO:03}. Lower bound]
The Mass Distribution Principle tells us that $\dimh E\geq s$, so $\dimh \clE_{\mathbf{t}}(B)\geq s$. Letting $s\to s_0(B)$, we conclude $\dimh \clE_{\mathbf{t}}(B)\geq s_0(B)$.
\end{proof}
\subsection{Proof of Theorem \ref{TEO:03}}\label{SUBSEC:PROOF:TEO:03}
\begin{proof}
The upper bound follows from Theorem \ref{LEM:TEO:Aux}. Certainly, for any $\veps>0$, every large $n\in\mathbb{N}$ satisfies $\Psi(n) \geq (B-\veps)^n$, so $\clE_{\mathbf{t}}(\Psi)\subseteq \clE_{\mathbf{t}}(B-\veps)$ and 
\[
\dimh\clE_{\mathbf{t}}(\Psi)
\leq
\dimh \clE_{\mathbf{t}}(B-\veps) \to \dimh \clE_{\mathbf{t}}(B) 
\;\text{ as }\; \veps\to 0.
\]
The lower bound is obtained in essentially the same way as in Theorem \ref{LEM:TEO:Aux}. The case $\frac{s_0(B)}{t_1} - \frac{2s_0(B)-1}{t_0}\leq 0$ is solved without significant modifications. Assume that
\begin{equation}\label{EQ:PF:TEO03:01}
\frac{s_0(B)}{t_1} - \frac{2s_0(B)-1}{t_0}>0.
\end{equation}
We shall only define a useful Cantor set contained in $\clE_{\mathbf{t}}(\Psi)$ and a probability measure supported on it. Let $\widetilde{B}>B$ be so close to $B$ that \eqref{EQ:PF:TEO03:01} still holds for $s_0(\widetilde{B})$. Let $\widetilde{A}$ be such that
\[
f_{t_0}\log \widetilde{A} = f_{t_0,t_1}(s)\log \widetilde{B}.
\]
Consider a strictly increasing sequence $(n_j)_{j\geq 1}$ in $\mathbb{N}$ such that
\[
\Psi(n_k)\leq \widetilde{B}^{n_k}
\quad\text{ for all }k\in\mathbb{N}.
\]
Write
\[
\beta_0\colon=\widetilde{A}^{1/t_0}
\quad\text{ and }\quad
\beta_1\colon=\widetilde{B}^{1/t_0}.
\]
Now, let $(\ell_k)_{k\geq 1}$ and $(i_k)_{k\geq 1}$ be the sequences of integers determined by
\[
(n_k-1) - (n_{k-1}-1)
=
\ell_kN + i_k 
\text{ with } 0\leq i_k< N.
\]
Call $\widetilde{E}$ be the subset of  $\clE_{\mathbf{t}}(\Psi)$ whose elements $x=\langle d_1, d_2,\ldots \rangle$ satisfy:
\begin{enumerate}[i.]
\item For each $k\in \mathbb{N}$, we have
\[
\beta_0^{n_k} 
\leq 
d_{n_k}
\leq  
2\beta_0^{n_k} 
\quad\text{ and }\quad
\beta_1^{n_k} 
\leq 
d_{n_k+1}
< 
2\beta_1^{n_k} . 
\]
\item We have $d_1=\ldots=d_{i_1}=2$ and $d_{n_k+2}=\ldots = d_{n_k + 1 + i_{k}}=2$ for $k\in \mathbb{N}$.
\item For every other natural number $n$, we have $2\leq d_n\leq M$.
\end{enumerate}
Hence, if $x=\langle d_1,d_2,\ldots, \rangle \in \widetilde{E}$, there are $\ell_k$ words $\bfw_1^{k}$, $\ldots$, $\bfw_{\ell_k}^{k}$ in $\{2,\ldots,m\}^N$ for $k\in\mathbb{N}$ such that
\[
\mathbf{d}
=
\underbrace{2\cdots 2}_{i_1 \text{ times }}\bfw_1^{1}\cdots \bfw_{\ell_1}^{1}d_{n_1}d_{n_1 +1}
\underbrace{2\cdots 2}_{i_2 \text{ times }}\bfw_1^{2}\cdots \bfw_{\ell_2}^{2}d_{n_2}d_{n_2 +1}
\ldots
\underbrace{2\cdots 2}_{i_k \text{ times }}\bfw_1^{k}\cdots \bfw_{\ell_k}^{k}d_{n_k}d_{n_k +1} 
\ldots .
\]
Define $\widetilde{D}\colon= \Lambda^{-1}\left[ \widetilde{E}\right]$. For each $n\in\mathbb{N}$ consider $\widetilde{D}_n \colon=\left\{ (d_1,\ldots, d_n) : (d_j)_{j\geq 1} \in \widetilde{D}\right\}$ and define the set $\widetilde{J}_n$ by adapting the definition of $J_n$ into our current context. 
Take any $1\leq n\leq n_1+1$ and $\mathbf{d}\in \widetilde{D}_n$. 
\begin{enumerate}[i.]
\item If $n\in\{1,\ldots, i_1\}$, put $\widetilde{\mu}(\widetilde{J}_n(2,\ldots,2))\colon=1$.
\item If $n=i_1 + 1 + N\ell$ with $1\leq \ell \leq \ell_1$, we write
\[
\widetilde{\mu}\left( J_n(\mathbf{d})\right) 
\colon=
\frac{1}{\beta_0^{sN\ell}}\left| I_{N\ell} (\bfw_1^1\cdots \bfw_{\ell}^1)
\right|^s.
\]
\item If $i_1+1 + N(\ell-1) +1 \leq n\leq i_1+1 + N(\ell-1) -1$, where $1\leq \ell \leq \ell_1$, write
\[
\widetilde{\mu} \left(J_n(\mathbf{d})\right)
\colon=
\sum_{\bfb} \widetilde{\mu}\left( J_{i_1 +1 +N\ell}(\mathbf{d}\bfb)\right),
\]
where the sum runs along all those words $\bfb$ such that $\mathbf{d}\bfb$ belongs to $\widetilde{D}_{i_1 +1 +N\ell}$.
\item When $n=n_1$, write 
\[
\widetilde{\mu}(J_{n_1}(\mathbf{d}))
\colon=
\frac{\widetilde{\mu}\left(J_{n_1-1}(d_1,\ldots, d_{n_1-1})\right)}{\lfloor 2\beta_0^{n_1}\rfloor - \lceil \beta_0^{n_1}\rceil}.
\]
\item When $n=n_1+1$, write
\[
\widetilde{\mu}(J_{n_1+1}(\mathbf{d})
\colon=
\frac{\widetilde{\mu}\left(J_{n_1-1}(d_1,\ldots, d_{n_1-1})\right)}{\left(\lfloor 2\beta_0^{n_1}\rfloor - \lceil \beta_0^{n_1}\rceil\right)\left(\lfloor 2\beta_1^{n_1}\rfloor - \lceil \beta_1^{n_1}\rceil\right)}.
\]
\end{enumerate}
Assume that we have already defined $\widetilde{\mu}$ for the fundamental intervals of order up to $n_k+1$ for some $k\in\mathbb{N}$. Take $n\in\mathbb{N}$ such that $n_k+2\leq n \leq n_{k+1}+1$.
\begin{enumerate}[i.]
\item If $n_k+2\leq  n\leq n_k+1 + i_{k}$, write
\[
\widetilde{\mu}\left( J_{n}(\mathbf{d})\right)
=
\widetilde{\mu}\left( J_{n_k+1}(d_1,\ldots, d_{n_k+1})\right).
\]
\item  If $n=n_{k} + i_{k+1} + 1 +  N\ell$ for some $1\leq \ell \leq \ell_k$, define
\[
\widetilde{\mu}\left( J_{n}(d_1,\ldots,  d_{n})\right)
\colon=
\frac{1}{2^{N\ell s}\beta_0^{sN\ell}}\,
\widetilde{\mu}\left( J_{n_k + i_{k+1} +1}(d_1,\ldots,  d_{n_k+1+i_{k+1}})\right).
\]
\item If $n_{k} + 1 +i_{k+1} + N(\ell-1)+1 \leq n < n_{k}+1+i_{k+1} + N \ell$, then
\[
\widetilde{\mu}(J_n(\mathbf{d}))
\colon=
\sum_{\bfb} \widetilde{\mu}\left(J_{n_{k}+1+i_{k+1} + N \ell}(\mathbf{d}\bfb)\right),
\]
where the sum runs along the words $\bfb$ such that $\mathbf{d}\bfb\in D_{n_{k}+1+i_{k+1} + N \ell}$.
\item If $n=n_k$, write
\[
\widetilde{\mu}(J_{n_1}(\mathbf{d}))
\colon=
\frac{\widetilde{\mu}\left(J_{n_k-1}(d_1,\ldots, d_{n_k-1})\right)}{\lfloor 2\beta_0^{n_k}\rfloor - \lceil \beta_0^{n_k}\rceil}.
\]
\item If $n=n_k+1$, write
\[
\widetilde{\mu}(J_{n_1+1}(\mathbf{d}))
\colon=
\frac{\widetilde{\mu}\left(J_{n_k-1}(d_1,\ldots, d_{n_k-1})\right)}{\left(\lfloor 2\beta_0^{n_k}\rfloor - \lceil \beta_0^{n_k}\rceil\right)\left(\lfloor 2\beta_1^{n_k}\rfloor - \lceil \beta_1^{n_k}\rceil\right)}.
\]
\end{enumerate}
\end{proof}

\section{Final remarks}\label{SEC:FINALREMARKS}

Our investigations give rise to a natural question: what happens when $1<B<\infty$ and $m\geq 3$? Unfortunately, our argument is not strong enough to solve this problem. However, based on \cite{BakHusKleWan2022}, we state a conjecture on the Hausdorff dimension of $\clE_{\mathbf{t}}(\Psi)$. For any $m\in\mathbb{N}$ and $\mathbf{t} = (t_0, \ldots, t_{m-1})\in \mathbb{R}_{>0}^m$, define the functions $f_{t_0}, f_{t_0, t_1}, \ldots,f_{t_0, t_1, \ldots, t_{m-1}}$ as follows: $f_{t_0}(s)=\frac{s}{t_0}$ and 
\[
f_{t_0,\ldots, t_{j}}(s)
= \frac{sf_{t_0,\ldots, t_{j-1}}(s)}{ t_jf_{t_0,\ldots, t_{j-1}}(s) + \max\left\{0,s - \frac{2s-1}{\max\{ t_0,\ldots, t_{j-1} \}  }\right\}}
\] 
for all $j\in\{2,\ldots, m\}$.
\begin{conj01}
Let $m\in\mathbb{N}_{\geq 3}$ be arbitrary and let $B$ be as in \eqref{EQ:DefBb}. If $1< B < \infty$, then $\dimh \clE_{\mathbf{t}}(\Psi)$ is the unique solution $s$ of
\[
\sum_{d=2}^{\infty} \frac{1}{d^s(d-1)^sB^{f_{t_0,\ldots, t_{m-1}}(s)}} 
 = 1.
\]
\end{conj01}
\bibliography{lur}
\bibliographystyle{abbrv}
\end{document}